\documentclass[12pt]{amsart}

\usepackage{amsmath,amsthm,amssymb,mathtools,url,graphicx,pdfpages,tikz,rotating, url}
\usepackage{float, comment, ytableau, enumitem, pgffor, circledsteps, etoolbox, xstring, mathtools, amsmath, amssymb, tikz, tikz-cd, adjustbox}
\usepackage[bookmarks=true,pdfborder={0 0 0}]{hyperref}
\usepackage[capitalize]{cleveref}
\usepackage[left=1in,right=1in,top=1.2in, bottom=1in]{geometry}
\newtheorem{theorem}{Theorem}[section]
\newtheorem{lemma}[theorem]{Lemma}
\newtheorem{corollary}[theorem]{Corollary}

\theoremstyle{definition}
\newtheorem{definition}[theorem]{Definition}
\newtheorem{notation}[theorem]{Notation}
\newtheorem{example}[theorem]{Example}

\newtheorem{remark}[theorem]{Remark}
\newtheorem{proposition}[theorem]{Proposition}

\numberwithin{equation}{section}

\newcommand{\y}[1]{\ydiagram{#1}}
\newcommand{\boks}[1]{\ytableausetup{boxsize = #1cm}}

\DeclareMathOperator{\sh}{\mathcal{Y}}
\DeclareMathOperator{\N}{\mathbb{Z}_{\geq 0}}
\newcommand{\la}{\lambda}

\DeclareMathOperator{\p}{\mathcal{P}}

\newcommand{\workhorse}{modular hook recursion}
\DeclareMathOperator{\g}{\mathcal{G}}

\DeclareMathOperator{\core}{{core}}
\DeclareMathOperator{\Vandermonde}{Van}
\renewcommand{\SS}{\mathfrak{S}}
\DeclareMathOperator{\specht}{Specht}

\DeclareMathOperator{\Od}{Od}
\ytableausetup{aligntableaux = center, nobaseline}
\allowdisplaybreaks

 \title{Enumeration of Odd-Dimensional Partitions Modulo 4}%write title

\author{Aditya Khanna}
\date{}
\begin{document}

%\author{Aditya Khanna}

{\let\newpage\relax\maketitle}

\begin{abstract}
The number of standard Young tableaux of a partition shape $\lambda$ is called the dimension of the partition $\la$. Partitions with odd dimensions were enumerated by McKay and were further characterized by Macdonald in the 1970s. Let $a_i(n)$ be the number of partitions of $n$ with dimension congruent to $i$ modulo 4. In this paper, we refine Macdonald's and McKay's results by computing $a_1(n)$ and $a_3(n)$ when $n$ has no consecutive 1s in its binary expansion or when the sum of binary digits of $n$ is 2.
\end{abstract}

\section{Introduction}\label{sec:introduction}
Representation theory offers a deep look into the inner workings of groups by studying their action on vector spaces. A vector space $V$ is called a \textit{representation} of $G$ if the elements of $G$ act like linear transformations on $V$. For a finite group, Maschke's theorem tells us that any complex representation decomposes into a direct sum of $G$-invariant subspaces called \textit{irreducible representations}, each indexed by a conjugacy class of $G$. The fact that the dimensions of irreducible representations are natural numbers leads one to question whether these numbers are counting some objects. An important family of objects of combinatorial consequence is \textit{partitions of $n$}. A partition of $n$ is a weakly decreasing list of natural numbers that sum to $n$.
Our focus in this paper is on the representation theory of $\SS_n$, the \textit{symmetric group on $n$ letters}, whose irreducible representations $\specht(\la)$ are indexed by partitions $\la$ of $n$. Every basis of $\specht(\la)$ is indexed by \textit{standard Young tableaux} (SYT) on the shape $\la$ (Prop. 2.5.9 and Thm. 2.6.4 in \cite{sagan}). The number of SYT of a given shape can be computed using the well-known hook-length formula \cite{frthook}.

The direction of inquiry that motivates this paper is the modular behavior of the dimensions (also known as \textit{degrees}) of the representations. Our story begins with the 2-page note of McKay \cite{mckay} wherein he, inspired by Brauer's assertion that ``there may be important properties of character tables yet to be discovered'', computes  the number of irreducible representations with odd dimensions for several groups. Let $m_p(G)$ for a prime $p$ and group $G$ denote the number of irreducible representations of $G$ whose dimensions are not divisible by $p$. The quantity $m_2(G)$ is often a power of 2 in McKay's computations, and takes a wonderfully concise form for the case of the symmetric group. If $n$ is a natural number such that $n = 2^{k_1}+\ldots + 2^{k_\ell}$ with $k_1 > k_2 >\ldots > k_\ell$, then $m_2(\SS_n) = 2^{k_1 + k_2 + \ldots + k_\ell}$. Thus, all of the information about the enumeration of representations of $\SS_n$ with odd dimensions is encoded in the binary expansion of $n$. The proof provided by McKay invokes the Murnaghan-Nakayama rule (\cite[Thm. 10.4.2]{sagan}), thereby reducing the problem to studying hooks of size a power of 2 on partition diagrams. In the same paper, McKay formulates  the celebrated \textit{McKay conjecture} which relates $m_2(G)$ for a group $G$ to $m_2(N)$ for the normalizer $N$ of a 2-Sylow subgroup of $G$. In his 1976 paper \cite{olsson-paper}, Olsson discusses conditions that imply McKay's conjecture for general primes $p$ and for $GL_q(n)$. This conjecture was proved in \cite{mckay-conjecture-odd} and a general form of the conjecture for all primes $p$ was established in a recent paper \cite{mckay-conjecture-full} by Cabanes and Sp\"ath. Navarro, an important contributor to the proof of the McKay conjecture, said that the decades-long work on this problem has resulted in ``beautiful, wonderful, deep mathematics.”  \cite{quanta} 

We are interested in the quantity $m_p(G)$ from a more enumerative and combinatorial perspective. The concise formula of McKay for $m_2(n)$ was generalized to an arbitrary prime $p$ through the use of \textit{$p$-core towers} by Macdonald \cite{mcd}. This proof method reduces the problem of finding the modularity of dimensions to the problem of labeling a $p$-ary tree with $p$-cores (\cref{def:core}). From this starting point, many avenues of extension present themselves. In \cite{mcd-coxeter}, Macdonald gave a formula for $m_p(G)$ when $G$ is a Coxeter group. In
\cite{jia-huang}, Huang extends the results for $m_p(\SS_n)$ to find the mod $p$ properties of projective indecomposable modules of a deformation of the symmetric group called the 0-Hecke algebra. All these enumeration results are for primes $p$. For a general $n$, asymptotic results \cite{peluse} show that $m_n(\SS_N) = 0$ as $N\to \infty$. 
%https://arxiv.org/pdf/2301.02203

Another direction of research, to which we contribute, tries to understand $m_n(\SS_N)$ explicitly for a natural number $n$. The 1976 paper \cite{olsson-paper} by Olsson studies heights of characters using Macdonald's theory of $p$-core towers which amounts to finding number of partitions with dimensions congruent to $p^i$ mod $p^n$ for $n\geq i$. This paper also contains positivity results for $m_2(G)$ (Prop. 3.9 in \cite{olsson-paper}) based on the relative size of the partition and its cores. From this proposition, we have an equivalent condition for the case $p = 2$, and thus an equivalent condition on the positivity of the number of partitions with dimensions 2 mod 4, but no explicit results. The explicit computations for the 2 mod 4 case using the 2-core tower method can be found in \cite{khanna-2mod4}. In this paper, we deal with the remaining cases modulo 4. 
Define $a_{i}(n)$  to be the number of partitions of $n$ with dimensions congruent to $i$ modulo $4$. Amrutha P and T. Geetha tackle the problem of computing $m_{2^k}(\SS_n)$ in their paper \cite{geetha}. They provide general recursive results and find $m_4(\SS_{2^\ell})$ and $m_8(\SS_{2^\ell})$, while also providing a characterization for partitions of $2^\ell$ with dimensions congruent to 2 modulo 4. 

The primary motivation for our project to compute $a_1(n)$ and $a_3(n)$ stems from the work \cite{jyo} of Spallone and Ganguly on the characterization of spinorial representations of $\SS_n$ and $\mathfrak{A}_n$ (the group of even permutations of $\{1,2,\ldots, n\}$). They prove that the spinoriality of $\specht(\la)$ is equivalent to the existence of a congruence condition modulo 8 involving the dimension of $\specht(
\la)$ and its value at a certain product of transpositions. As mentioned in the introduction to their paper, spinoriality of Galois representations has applications to number theory.

In this paper, we establish results for $a_1(n)$ and $a_3(n)$ for sparse numbers (  \cref{def:sparse}), which when combined with results for $a_2(n)$ from \cite{khanna-2mod4} allows us to compute $m_4(
\SS_n)$ for all sparse numbers. Furthermore, as outlined in \cite{khanna}, these results can be used to find $m_2(\mathfrak{A}_n)$. Our attempt to solve the mod 4 problem for partition dimensions is a foray into understanding the modulo 8 properties of the dimension, and thus better understanding the spinoriality of symmetric and alternating group representations. This problem is also of independent interest as it showcases the difficulty inherent in computing $m_n(G)$ for general $n$. We only capture a sliver of this expansive enumeration endeavor, but we do so with a reinforcement of the same assertion that inspired McKay to investigate this field in the first place --  ``\textit{there are many properties of character tables yet to be discovered}''.

\subsection{Background} 
We split our discussion of the necessary concepts and definitions for this paper into two parts. The first is this section where we discuss the well-known combinatorial concepts that are needed to understand the main problem tackled in the paper. In \cref{sec:notation}, we present the definitions necessary to understand the main proof method. We start with the definition of a partition.

We call a $k$-tuple of positive integers $\lambda = (\lambda_1, \ldots, \lambda_k)$ a \textit{partition} of \textit{size} $n$ if $\lambda_1 + \lambda_2 + \ldots + \lambda_k= n$. We denote by $\lambda\vdash n$ that $\la$ is a partition of $n$. We denote the size of a partition $\la$ by $|\la| = n$ and the \textit{length of the partition} by $\ell(\lambda) = k$.
We construct the \textit{Ferrers diagram of $\lambda$}, denoted by $\sh(\lambda)$, by placing $\lambda_i$ top-left justified boxes in the $i^{\text{th}}$ row from top.
For  $\lambda = (3, 3, 2, 1)$, we have \boks{0.3}
\[
\sh(\lambda) = \ydiagram{3, 3, 2, 1}
\]
\ytableausetup{boxsize = normal}
For $\lambda\vdash n$, we can fill the boxes of $\sh(\lambda)$ with numbers in the set $\{1, 2, \ldots, n\}$. A \textit{standard Young tableau (SYT)} on the shape $\lambda$ is a filling of $\sh(\la)$ such that the entries increase strictly along rows (left to right) and along columns (top to bottom).
Continuing with the above example, the following filling is an SYT on $(3,3,2,1)$:
\boks{0.42}
\ytableausetup{nobaseline}
\begin{center}
\ytableaushort{1 2 7, 3 5 8, 4 9, 6}
\end{center}

The number of SYTs of shape $\lambda$ is denoted by $f^{\lambda}$. We call this quantity the \textit{dimension of $\la$}, and in the literature the terminology \textit{degree of $\la$} also appears \cite{mckay}. We can explicitly compute $f^\la$ using the \textit{hook-length formula} (\cref{eq:frt}) which requires us to first define hook-lengths.
We call the boxes in Ferrers diagrams \textit{cells}. We label the cell in the $i^{\mathrm{th}}$ row from top and $j^{\text{th}}$ column from the left by $(i,j)$. Write $(i,j) \in \sh(\lambda)$ if $1\leq i\leq \ell(\lambda)$ and $1\leq j \leq \lambda_i$. Define the \textit{border of $\lambda$} to be the set of all cells $(i,j)\in \sh(\lambda)$ such that at least one of $(i, j+1)$, $(i+1, j)$ or $(i+1, j+1)$ is not contained in $\sh(\lambda)$. 
\boks{0.3} The partition $(5,5,5,4,2)$ has the border shaded in gray as below:
    \begin{center}
\ydiagram{5,5,5,4,2}*[*(lightgray)]{4+1,4+1,3+2,1+3,2}
\end{center}

For a cell $(i,j)$ in $\sh(\lambda)$, define the \textit{rim hook
%\footnote{In the literature, the objects we define are known as \textit{rim hooks} while the terminology \textit{hooks} is reserved for a related concept. The sets of hooks and rim hooks are in bijection, so we use the term ``hooks" in place of ``rim hooks" for convenience.}
of $\lambda$ at $(i,j)$} to be the set of border cells weakly below row $i$, and weakly to the right of column $j$. The number of cells contained in the rim hook at $(i,j)$ is called its \textit{hook-length} and is denoted by $h_{i,j}$. A rim hook of length $t$ is called a $t$-rim hook. The partition $(5,5,5,4,2)$ has a $5$-rim hook (shaded in gray) at $(2,3)$ (shaded in black):
\begin{center}
\ydiagram{5,5,5,4,2}*[*(lightgray)]{0, 4+1, 3+2, 2 + 2}*[*(black)]{0,2+1}
\end{center} 
It is through these hook-lengths that the computation of $f^\la$ is made possible. From \cite{frthook}, we know
\begin{equation}\label{eq:frt}
    f^\lambda = \frac{n!}{\prod\limits_{(i,j)\in \sh(\lambda)} h_{i,j}}.
\end{equation}
In \cref{hlfsec}, we will describe an equivalent formula due to Frobenius which will be more amenable to the methods of this paper. Our methods will require us to \textit{remove} the rim hooks from partitions and see how the mod 4 behavior of $f^\la$ changes.
As the process of rim hook removal and the output of this process are so crucial, we describe them now in terms of Ferrers diagrams, and later in terms of sets in \cref{hookremprop}. 
In terms of diagrams, we \textit{remove a $t$-rim hook $H$ from $\lambda$} by removing all cells of $H$ from $\sh(\lambda)$.
Removing the 5-rim hook at $(2,3)$ from $(5,5,5,4,2)$ results in the partition $(5,4,3,2,2)$ with the diagram \boks{0.3}
\begin{center}
\ydiagram{5,4,3,2,2}
\end{center} 

If the removal of a $t$-rim hook $H$ from $\la$ results in the partition $\mu$, then we write $H = \la/\mu$.
\begin{definition}[$t$-core]\label{def:core}
If $\sh(\mu)$ does not contain any $t$-rim hook, then we call $\mu$ a \textit{$t$-core}.
\end{definition}

Given a partition $\lambda$, successive removals of $t$-rim hooks result in a partition that is eventually a $t$-core. The $t$-core thus obtained is independent of the order in which $t$-rim hooks are removed \cite[Thm. 2.7.16]{james-kerber}, and is called the \textit{$t$-core of $\lambda$}. We denote it by $\core_t(\lambda)$. The main proof method of this paper focuses on the mod 4 relationship between $f^\la$ and $f^\mu$ where $\mu$ is the $2^R$-core of $\la$, and $2^R$ is the largest power of 2 smaller than the size of $\la$.

%We provide explicit enumeration results in this paper. We present our results now and data for our quantity of interest, $\delta$, can be found \href{https://drive.google.com/file/d/1iaQz0ffXwzNc79LFg1_fDRkH-BLWTtL7/view?usp=drive_link}{here}.

\subsection{Main results}\label{ssec:mainresults}
We denote the set of partitions of $n$ with dimensions congruent to $i$ mod 4 by $a_i(n)$ throughout this paper. The number of odd-dimensional partitions is denoted by $a(n)$, which we may compute as $a(n) = a_1(n) + a_3(n)$. By the work of McKay and Macdonald, we know that $a(n) = 2^{k_1 + k_2 + \ldots + k_l}$ where $2^{k_i}$ appear in the binary expansion of $n$. In the notation of McKay described in the introduction, $a(n)= m_2(\SS_n)$. Furthermore, $m_4(n) = a_1(n) + a_2(n) + a_3(n)$. The problem of computing $a_2(n)$ is tackled in \cite{khanna-2mod4} while in this paper we aim to find $a_1(n)$ and $a_3(n)$.
The computation of $a_1(n)$ or $a_3(n)$ is a hard question in general, and thus, we prove the results for sparse\footnote{These numbers are called \textit{fibbinary numbers} in the literature as they are enumerated by the Fibonacci numbers. The reader may refer to the OEIS entry \texttt{A003714} \cite{oeis}.} numbers, which are positive integers with no consecutive ones in their binary expansions. For instance, $41$ is sparse as its binary expansion is $\mathtt{101001}$.

The methods described in this paper compute the difference $\delta(n) = a_1(n) - a_3(n)$ and the formulas for $\delta(n)$ are more concise than the formulas for $a_1(n)$ and $a_3(n)$ themselves. We can recover $a_1(n)$ and $a_3(n)$ from $\delta(n)$ and $a(n)$ using $a_1(n) = \frac{a(n) + \delta(n)}{2}$ and $a_3(n) = \frac{a(n)-\delta(n)}{2}$, respectively. We now present the recursions and explicit computations, which are proved in this paper.
\begin{theorem}\label{mainthm}
Let $n, m, R \in \mathbb{N}$ with $R\geq 2$ and $m > 0$. Suppose, $n = 2^{R} + m$ with $m<2^{R-1}$. Then we have
\[
\delta(n) = \begin{cases}
0, & \text{if }n \text{ is even}\\
4 \delta(m), &\text{if } n \text{ is odd}.
\end{cases}
\]
\end{theorem}

Let $\nu(n)$ denote the number of 1s in the binary expansion of $n$. For instance, $\nu(41) = 3$.
\begin{corollary}\label{maincor}
If $n$ is a sparse number, then
\[
\delta(n) = \begin{cases}
2\text{,} & \text{if }n = 2\\
0, &\text{if } n > 2 \text{ is even}\\
4^{\nu(n) - 1}, & \text{if } n\text{ is odd},
\end{cases}
\]
\end{corollary}
In the case where the binary expansion starts with two $\mathtt{1}$s followed by all zeroes, i.e., $n =\texttt{1100}\ldots$, we have the following result.
\begin{theorem}\label{11thm}
Let $n = 2^R + 2^{R-1}$ with $R\geq 1$. Then 
\[
\delta(n) = \begin{cases}
2, & \text{if } R = 1\\
8, & \text{if } R = 2\\
0, & \text{otherwise}.
\end{cases}
\]
\end{theorem}
We prove \cref{mainthm} and \cref{maincor} in \cref{sec:proofmain}, and \cref{11thm} in \cref{sec:proof11}. Before we get into the technical definitions and proofs, we give a broad overview of our proof ideas.

%%%%%%%%%%%%%%%%%%%%%%%%%%%%%%%%%%%%%%%%%%%%%%%%%%
 \subsection{Proof strategy and outline}\label{ssec:outline}
Our final goal in this paper is to compute $\delta(n)$ for non-negative integers $n$. Say we have a function $\Od$ on natural numbers which outputs $+1$ when the input is $1$ mod 4 and $-1$ when the input is $3$ mod 4. This allows us to use $\Od(f^\la)$ to associate a sign $+1$ with each partition with dimension 1 mod 4, and $-1$ with each partition with dimension 3 mod 4. For a fixed $n$, this process results in $a_1(n)$ partitions of $n$ with associated sign $+1$ and $a_3(n)$ partitions of $n$ with the sign $-1$. This allows us to write $\delta(n) = a_1(n) - a_3(n)$ as $\sum_{\la} \Od(f^\la)$ where the sum is over all odd-dimensional partitions $\la$. The technique we adopt to compute this sum is recursion arising from the process of rim hook removal (\cref{sec:workhorse}). If $n = 2^R + m$ with $m < 2^R$, then Macdonald's work \cite{mcd} tells us that each odd-dimensional partition $\la$ of $n$ yields an odd-dimensional partition $\mu$ of $m$ the removal of a unique $2^R$-rim hook, or in other words, $\mu$ is the $2^R$-core of $\la$. With this process of $2^R$-rim hook removal, we can associate a sign of $\pm 1$ computed by the ratio $\frac{\Od(f^\la)}{\Od(f^\mu)}$. Furthermore, define $S(\mu)$ as the sum $\sum_\la \frac{\Od(f^\la)}{\Od(f^\mu)}$ over all odd-dimensional partitions $\la$ of $n$ that yield $\mu$ after the $2^R$-rim hook removal process. When $\Od(f^\la)/\Od(f^\mu)$ is $+1$, then the dimension of $\la$ has the same residue modulo 4 as the dimension of $\mu$, and the residue is different when the ratio is $-1$. Thus, the quantity $S(\mu)$ captures the difference between the number of partitions $\la$ whose dimensions have the same or different residue modulo 4 compared to $f^\mu$. This allows us to write $\delta(n) = \sum_{\la} \Od(f^\la)$ with the sum over all odd-dimensional partitions $\la$ of $n$ as $\sum_{\mu} S(\mu) \Od(f^\mu)$ where the sum is now over all odd-dimensional partitions of $m$. 
The bulk of our paper will be devoted to understanding $S(\mu)$ and $\frac{\Od(f^\la)}{\Od(f^\mu)}$ for odd-dimensional partitions $\mu$ of $m$ and $\la$ of $n$ as above. For a general $n$, this problem is difficult, as the value of $S(\mu)$ can depend on the parts of $\mu$ (\cref{sec:problems}). In this paper, we opt to study the special cases where $n = 2^R + m$ with $m < 2^{R-1}$(\cref{sec:proofmain}), and the case where $n = 2^R + 2^{R-1}$ for $R\geq 1$ (\cref{sec:proof11}). The usefulness of these specializations becomes evident in the fact that $S(\mu)$ for $\mu \vdash m$ or $\mu \vdash 2^{R-1}$ is independent of the partition $\mu$ and only depends on the parity of $m$ (\cref{thm:parent-sum}). In these cases, we may factor out some constant $s_m$ dependent only on $m$ from the sum, allowing us to construct a very simple recursion for $\delta(n)$ which is $\delta(n) =s_m\sum_{\mu} \Od(f^\mu)=s_m\delta(m)$ where the sum is over all odd-dimensional partition $\mu$ of $m$.

The implementation of the above proof strategy will require technical arguments relating to objects arising from partitions called $\beta$-sets. In \cref{sec:notation}, we introduce these objects, their effectiveness in describing rim hook removal, and the Frobenius hook-length formula that uses $\beta$-sets to compute the dimension of partitions. We also present the characterization of odd-dimensional partitions due to Macdonald \cite{mcd} in terms of $\beta$-sets, and introduce the notion of $t$-\textit{parents}. This section is also where we properly introduce the function $\Od$ and lay out some of its required properties. In \cref{sec:workhorse}, we prove the \workhorse\, that establishes $\Od(f^\la)/\Od(f^\mu)$ for partitions $\la$ of $n = 2^R + m$ and $\mu$ of $m$.  In \cref{sec:proofmain}, we focus on the case of sparse numbers, and establish the recursion relating $\delta(n)$ and $\delta(m)$ (\cref{mainthm}) when $m < 2^{R-1}$. In \cref{sec:proof11}, we tackle the simplest non-sparse case where $n = 2^R + 2^{R-1}$ to prove \cref{11thm}. The proof in this section resembles the proofs presented in \cref{sec:proofmain}, but as we shall see, the casework required in this case is indicative of the general difficulty of the modulo 4 problem as seen through the lens of our methods. In \cref{sec:general-n}, we provide the interested reader with some pointers that might be helpful in resolving the general cases, and the obstructions that one may face. In \cref{sec:problems}, we provide data for $\delta(n)$ when $n < 128$, and comment on some patterns within it. At the beginning of each subsection, we provide some exposition to motivate our methods and situate the reader as the proofs are technical and computationally involved. 

We now summarize the specializations considered in this paper.
Throughout, we write $n = 2^R + m$ for $0 \leq m < 2^R$, where $2^R$ is the largest power of 2 smaller than $n$, and $m$ is what remains after subtracting this power of 2 from $n$.
In \cref{sec:notation} and \cref{sec:workhorse}, we consider the general case of $m < 2^R$ as all our propositions hold without any restrictions. In \cref{sec:proofmain}, we specialize to the case $m < 2^{R-1}$ which allows for explicit recursion formulas and enumeration. In \cref{sec:proof11}, we take the first step beyond the non-sparse case by considering $m = 2^{R-1}$. In \cref{sec:problems}, we discuss strategies for the case where $2^{R-1} < m < 2^{R}$.
 
\section{Development of the modular hook recursion formula}\label{sec:notation}
The main recursion underlying our methods requires some background to set up. In this section, we present an overview of the foundation needed to derive the formula in the next section and understand the proof methods in the rest of the paper. The upcoming definitions hold for all values of $n = 2^R + m$, $m < 2^R$, and we will in later sections impose further restrictions on $m$.

Let $n = 2^R + m$ for natural numbers $n$, $R$ and $m$ with $m < 2^R$.
To compute $\delta(n)$, we study the modulo 4 properties of the Frobenius hook-length formula (\cref{HLF}) which counts the number of SYT $f^\la$ of a given shape $\la$. The main ingredient of this formula is the set of hook-lengths which appear in the cells of the first column of the diagram (\cref{notn:first-column}). For an odd-dimensional partition $\la$, we keep track of the residue modulo 4 of $f^\la$ through the function $\Od(f^\la)$ (\cref{odsubs}) which outputs $+1$ when $f^\la$ has residue 1 mod 4 and $-1$ when $f^\la$ has residue 3 mod 4.
To establish the main recursion in \cref{sec:workhorse}, we explore how the removal of a $2^R$-rim hook from a diagram $\sh(\lambda)$ affects the value of $\Od(f^\lambda)$. The process of hook-removal can be understood naturally in terms of $\beta$-sets, which motivates our usage of them as the central object in our paper.

\subsection{\texorpdfstring{$\beta$}{beta}-sets}
The Frobenius formula described in the next section takes as input the first column hook-lengths of a partition and outputs its dimension. The formula works for a general class of equivalent objects called $\beta$-sets which we now describe. These $\beta$-sets are a natural candidate for our purposes as removal of a rim hook on diagrams can be understood as a simple alteration of a single element in a $\beta$-set as described in \cref{hookremprop}.
\begin{notation}\label{notn:first-column}
    The set of first column hook-lengths of $\lambda$ is denoted by $H(\lambda) := \{h_{i, 1}\mid 1\leq i \leq \ell(\lambda)\}$,
where we conventionally list the elements in descending order.  Explicitly, if $\lambda = (\lambda_1, \ldots, \lambda_k)$, then $h_{i,1} = \lambda_i + k - i$. In this paper, we often drop the column subscript 1 and just write $h_i$ for $h_{i,1}$. 
\end{notation}
For any finite $X\subset \N$, define the \textit{$r$-shift of $X$} to be \[X^{+r} = \{x + r\mid x\in X\} \cup \{0, \ldots, r-1\}.\] Also, define $X^{+0} = X$. Thus, an $r$-shift of $X$ increases each element of $X$ by $r$ and appends the elements $0,1,\ldots, r-1$ to this new set. For a partition $\lambda$, sets of the form $H(\lambda)^{+r}$ for $r\geq 0$ are known as the \textit{$\beta$-sets of $\lambda$}. We now present some important properties of $\beta$-sets as remarks.

\begin{remark}
    For \textit{any} weakly decreasing $r$-tuple of non-negative integers $w$, we define $H(w) = \{w_i + r - i : 1\leq i\leq r\}$.
    An alternate way to view the $\beta$-set $H(\lambda)^{+r}$ is by using partitions with trailing zeroes. If $r\geq \ell(\lambda)$, then let $\rho$ be the weakly decreasing list of integers obtained by appending $r-\ell(\lambda)$ zeroes to the end of $\la$. If we compute $H(\rho) = \{h_1, \ldots, h_r\}$ using the formula $h_i = \rho_i + r - i$, we obtain exactly $H(\lambda)^{+(r-\ell(\la))}$. For instance, let $\la = (4,3,3,1)$ with $H(\la) = \{7,5,4,1\}$. If $r = 6$, then $H(\lambda)^{+2} = \{9,7,6,3,1,0\}$. On the other hand, we have $\rho = (4,3,3,1,0,0)$ and we find that $H(\rho)$ agrees with $H(\la)^{+2}$.
\end{remark}

\begin{remark}\label{max-of-H-lambda}
    As the largest hook-length can at most encompass all the cells of a partition $\la$, the maximum element of $H(\lambda)$ is at most $|\lambda|$. In general, the maximum element of an $r$-shift of $\la$ is at most $|\la| + r$.
\end{remark}
\begin{remark}\label{rem:partition-from-beta-set}
     We can reconstruct the original partition from a $\beta$-set with $k$ elements by setting $\lambda_i = h_i + i - k$ and considering only the values where $\lambda_i >0$. For instance, if we have a $\beta$-set $X = \{8,4,3,1,0\}$ of length $k= 5$, then the corresponding partition $\la$ has non-negative parts $\la_1 = 8 - 5 + 1 = 4$, $\la_2 = 4 - 5 + 2 = 1$, $\la_3 = 3 - 5 + 3 = 1$, $\la_4 = 1 -5 + 4 =0$ and $\la_5 = 0 -5 + 5 = 0$. Thus, $\la = (4,1,1)$ and $X$ is $H(\la)^{+2}$. 
\end{remark}
From the above remark, we see that if two sets can be expressed as $r$-shifts of each other, then they must arise as first column hook-length sets of the same partition. This motivates the following equivalence relation.
\begin{notation}\label{notn:sim}
    For finite sets $X, Y\subset \mathbb{N}$, we say $X\sim_{\beta} Y$ if and only if $X = Y^{+r}$ or $Y = X^{+r}$ for some $r \in \mathbb{Z}_{\geq 0}$. 
\end{notation}

We now explicitly describe how to interpret $t$-rim hook removal and addition in terms of {$\beta$-sets}. Given a set $S$ with $a\in S$ and $b\not\in S$, define $S[a\to b]$ to be the set obtained from $S$ where $a$ is replaced by $b$. In notation, $S[a\to b] = (S\cup \{b\})\backslash\{a\}$. In the cases where $b$ is already in $S$ or $a$ is not in $S$, we define $S[a\to b] = S$. 
\begin{lemma}\label{lem:double-replacement}
    For a set $S\subset \mathbb{Z}_{\geq 0}$ with $a\in S$ and $b\in \mathbb{Z}_{\geq 0}$ which may or may not be in $S$, we have
$S[a\to b][b\to a] = S$.
\end{lemma}
\begin{proof}
    If $b\not\in S$, then $S[a\to b]$ contains all the same elements as $S$ except with $a$ replaced by $b$. The set $S[a\to b][b\to a]$ replaces the $b$ back by $a$, giving us $S$. On the other hand, if $b\in S$, then $S[a\to b] = S$ and $S[b\to a] = S$.
\end{proof}

\begin{proposition}\label{hookremprop}
Let $\lambda$ be a partition and $X$ be a $\beta$-set of $\lambda$. The diagram $\sh(\lambda)$ contains a $t$-rim hook if and only if there exists an $h$ in $X$ such that $h \geq t$ and $h - t$ is not in $X$. Furthermore, if $\mu$ is the partition obtained after removing the $t$-rim hook, then $H(\mu) \sim_{\beta} H(\lambda)[h\to h-t]$.
\end{proposition}
\begin{proof}
See Corollary 1.5 on page 7 of \cite{olsson} or Lemma 2.7.1 of \cite{james-kerber}.
\end{proof}
\cref{hookremprop} states that the removal of a $t$-rim hook from the diagram of $\la$ is expressed in terms of the $\beta$-set $X$ of $\la$ by selecting an element $h\in X$ and reducing it by $t$ (if $h-t$ is not already an element of $X$). Equivalently, one may add a $t$-rim hook to $\sh(\la)$ by selecting an element $h'$ in a sufficiently large $\beta$-set $Y$ of $\la$ and increasing it by $t$ (if $h'+t$ is not in $X$). In the case where no value $h$ in $X$ can be reduced by $t$ as $h-t$ already occurs in $X$,  we declare the partition $\la$ to be a $t$-core.
\begin{example}
In this example, we find the 5-core of $\lambda = (6,5,5,4,2)$ using the process of rim hook removal described on $\beta$-sets, and then display the same process using diagrams. We begin with $H(\lambda) = \{10, 8, 7, 5, 2\}$. In this case, $8 \in H(\lambda)$ but $3 \not\in H(\lambda)$. Thus, there exists a 5-rim hook in $\sh(\lambda)$ which we can remove to obtain $\lambda'$ (not to be confused with the transpose) such that $H(\lambda') = \{10, 7, 5, 3, 2\}$. Now, we can further replace $5$ by 0 to obtain $H(\lambda'') = \{10, 7, 3, 2, 0\} \sim_{\beta} \{9,6,2,1\}$. Finally, we can replace 9 by 4 to get $H(\lambda''') = \{6,4,2,1\}$.
No further reduction of elements by 5 can be performed. Thus, $\lambda''' = (3,2,1,1)$ is the 5-core of $\lambda$. This can be visualized using Ferrers diagrams with the 5-rim hooks removed at each step shaded in gray:
\[
\lambda = \y{6,5,5,4,2}*[*(lightgray)]{0,4+1,3+2,2+2} \to \lambda' = \y{6,4,3,2,2}*[*(lightgray)]{0,0,1+2,1+1,2} \to \lambda'' = \y{6,4,1,1}*[*(lightgray)]{3+3,2+2} \to \lambda''' = \y{3,2,1,1}
\]
\end{example}

\subsection{Characterization of odd-dimensional partitions}\label{hlfsec}
The following proposition is due to Frobenius:
\begin{proposition}[Eq. 6 in \cite{frobenius}]\label{HLF}
Let $\lambda$ be a partition of $n$ and let $X = \{x_1, \ldots, x_k\}$  with $x_1 > \ldots > x_k$ be any $\beta$-set of $\lambda$. Then
\[
f^{\lambda} = \dfrac{n!\prod\limits_{1\leq i<j\leq k} (x_i - x_j)}{\prod\limits_{i=1}^k x_i!}.
\]
\end{proposition}
Note that in the above proposition $X$ is any $\beta$-set, and not necessarily the set of first-column hook-lengths, and can possibly contain a 0.
\begin{remark}\label{rem:three-parts}
    The above expression for $f^\la$ can be split into three parts: the factorial of the size of the partition, $n!$; the product of the factorials of the elements of the $\beta$-set, $\prod_{i} x_i!$; and the Vandermonde product, $\Vandermonde(x_1, x_2, \ldots, x_k) :=  \prod\limits_{1\leq i<j\leq k} (x_i - x_j)$.
\end{remark}

\begin{example}\label{ex:f-la}
    Let $\lambda = (4,4,2,1,1)$ with $\beta$-set $X=\{10,9,6,4,3,1,0\}$. Then \[f^\lambda = \dfrac{12! \cdot (10-9)(10-6)(10-4)\ldots (3-1)(3-0)(1-0)}{10!\,9!\,6!\,4!\,3!\,1!\,0!}=4455.\]
\end{example}

We call $\lambda$ an \textit{odd-dimensional partition} if $f^{\lambda}$ is odd. For instance, the partition $(4,4,2,1,1)$ in \cref{ex:f-la} is an odd-dimensional partition.
We now restate the characterization of odd-dimensional partitions by Macdonald~\cite{mcd} in terms of rim hooks and cores. 

\begin{proposition}\label{mcdprop}
If $\lambda$ is a partition of $n = 2^R + m$ with $m < 2^R$, then
$\lambda$ is an odd-dimensional partition if and only if $\lambda$ contains exactly one $2^R$-rim hook and $\core_{2^R}(\lambda)$ is also an odd-dimensional partition.
\end{proposition}
\begin{proof}
This proposition is stated as Lemma 1 in \cite{aps} and a proof of it is given in Section 6 of the same paper.
\end{proof}

\begin{remark}\label{odd-partitions-are-hooks}
    A consequence of the above proposition is that $\la$ is an odd-dimensional partition of $2^R$ if and only if $\la$ is a $2^R$-rim hook.
\end{remark}
The repeated application of \cref{mcdprop} tells us that if we can remove a unique rim hook of size the largest power of 2 smaller than the size of the partition at each step, then the partition is odd-dimensional. For instance, if $\la$ is a partition of $2^{6} + 2^{3} + 2^1$, we can remove a unique $2^6$-rim hook from $\la$ to yield $\mu$, a unique $2^3$-rim hook from $\mu$ to yield $\nu$, and a unique 2-rim hook from $\nu$ (which is $\nu$ itself by the above remark) to get $\varnothing$. As might be evident from this example, binary expansions play an important role in our techniques, so we present the notation for them that is used throughout the text.
\begin{notation}[Binary expansion]\label{not:binary} Let $n = \sum\limits_{i=0}^k b_i 2^i$ with $b_i \in \{0,1\}$ such that $b_k = 1$. We call the $b_i$  the \textit{bits of $n$}. We use the shorthand $n = [b_k\ldots b_0]_2$ in place of the expansion.
\end{notation}
Each odd-dimensional partition can be constructed by repeated addition of larger rim hooks corresponding to each positive bit in its binary expansion starting from the empty shape $\varnothing$. In the above example, we removed rim hooks to get to the empty partition, but in our proofs, we start from a partition, and then add rim hooks in all possible ways. We call the partitions obtained by addition of $t$-rim hooks to $\sh(\mu)$ the $t$-\textit{parents} of $\mu$.
\begin{definition}\label{parentsrem}
Let $\lambda$ be an odd-dimensional partition of $n = 2^R + m$ with $m < 2^R$. If $\core_{2^R}(\lambda) = \mu$, then we call $\lambda$ a \textit{$2^R$-parent of $\mu$}. 
\end{definition}
\cref{parentsrem} is unambiguous as $\lambda$ contains a unique $2^R$-rim hook which when removed results in a partition of size $m < 2^R$ which is a $2^R$-core. Equivalently, $\la$ is a $2^R$-parent of $\mu$ if $\la$ is obtained by the addition of a single $2^R$-rim hook to $\mu$.  Recall from \cref{hookremprop} that the addition (resp. the removal) of a rim hook to $\mu$ corresponds to increasing (resp. decreasing) an element in the $\beta$-set of $\mu$. We can create $\beta$-sets of $2^R$-parents of an odd-dimensional partition $\mu$ by starting with a $\beta$-set $B$ of $\mu$ and replacing some element $b$ in $B$ by $b+2^R$. 

Fix a partition $\mu$. The $2^R$-parents of $\mu$ can be split into two categories: Type I $2^R$-parents are formed by choosing the $\beta$-set $H(\mu)$ and replacing some $h\in H(\mu)$ by $h+2^R$, whereas Type II $2^R$-parents arise from replacing $0$ by $2^R$ in $\beta$-sets which are $r$-shifts of $H(\mu)$ for $1\leq r \leq 2^R$. To understand how the $\beta$-sets of $2^R$-parents of $\mu$ can be derived explicitly from the $\beta$-sets of $\mu$, we present the next proposition.
\begin{proposition}\label{parents}
If $\mu \vdash m < 2^R$ is obtained from $\lambda\vdash m+2^R$ after removal of a $2^R$-rim hook, then exactly one of the following holds:
\begin{enumerate}
\item (Type I) $H(\lambda) = H(\mu)[x\to x+2^R]$ for some $x\in H(\mu)$,
\item (Type II) $H(\lambda) = H(\mu)^{+r}[0\to 2^R]$ for $1\leq r\leq 2^R$ such that $2^R\not\in H(\mu)^{+r}$.
\end{enumerate}
\end{proposition}
\begin{proof}
By \cref{hookremprop} and \cref{mcdprop}, recall that the $\beta$-set obtained after a removal of a $2^R$-rim hook from $H(\la)$ results in a $\beta$-set of $\mu$. In the notation from \cref{notn:sim},  we have $H(\lambda)[h \to h-2^R]\sim_{\beta}H(\mu)$. As $\lambda$ is formed by addition of a $2^R$-rim hook to $\mu$, $\lambda$ contains at least as many parts as $\mu$, which means $\ell(\la) \geq \ell(\mu)$. This turns the equivalence $H(\lambda)[h \to h-2^R]\sim_{\beta}H(\mu)$ into the equality $H(\lambda)[h\to h-2^R] = H(\mu)^{+r}$ where $r\geq 0$. As explained in \cref{rem:typeii-quirks} below, the $r$ takes values in $\{0,1,2,\ldots,2^{R}\}$. When $r =0$, we get $H(\lambda)[h\to h-2^R] = H(\mu)$ where $h$ is an element of $H(\la)$ and $h-2^{R}$ is an element of $H(\mu)$. If we let $h - 2^R = x$, we obtain the equality $H(\la)[x+2^R \to x] = H(\mu)$ for $x\in H(\mu)$. We also have
\[
H(\la)[x+2^R \to x][x\to x+2^R] = H(\mu)[x\to x+2^R]
\] 
which gives $H(\la)= H(\mu)[x\to x+2^R]$ \cref{lem:double-replacement}. 
Following \cref{rem:typeii-quirks} below, we may take $h - 2^R = 0$ without loss of generality. Thus, we have for some $0 < r \leq 2^R$, $H(\la)[2^R\to 0] = H(\mu)^{+r}$, which by another application of \cref{lem:double-replacement} gives the desired result. We emphasize that as $2^R$ is being replaced in $H(\la)$ by 0, $2^R$ cannot be an element in $H(\mu)^{+r}$, thus we get no corresponding $2^R$-parents for the $r_h$-shift where $r_h = 2^R-h$ for some $h\in H(\mu)$, leading to the condition in part 2.
\end{proof}
\begin{remark}\label{rem:typeii-quirks}
   In \cref{parents}, $r$ takes values only in the from 1 to $2^R$ because any $\beta$-set $Y = H(\mu)^{+r}$ with $r > 2^R$ already contains the element $2^R$, thus $Y[0\to 2^R] = Y$ would not result in a $2^R$-parent.  
    In the case of Type II $2^R$-parents, instead of replacing 0 by $2^R$, we may consider the general case of replacing any $h\in H(\mu)^{+r}$ with $h+2^R$. It follows from the definition of $r$-shifts that the $\beta$-set $X = H(\mu)^{+r}[h\to h+2^R]$ with $h < r$ is equal to $X = (H(\mu)^{r-h}[0\to 2^R])^{+h}$ which makes $X$ equivalent to $H(\mu)^{r-h}[0\to 2^R]$. Thus the $2^R$-parent of the $\beta$-set $X$ is the same as that of $H(\mu)^{r-h}[0\to 2^R]$. Only replacing the element 0 by $2^R$ in the shifted $\beta$-sets ensures that each Type II $2^R$-parent is counted exactly once.
\end{remark}
\begin{remark}
    In \cite{inv_kostka}, the Type II $t$-parents arise from the addition of \textit{special} $t$-rim hooks, which are $t$-rim hooks whose southwestern-most cell lies in the first (leftmost) column.
\end{remark}

%Feb 10
\begin{example}
We use the notation from \cref{parents}.
Let $\mu = (4,2)$ be an odd-dimensional partition and let $2^R = 8$. We have $H(\mu) = \{5,2\}$. We construct the Type I 8-parent $\la^{(2)} = (9,5)$ corresponding to $x= 2$ by adding 8 to the entry $2\in H(\mu)$ and obtaining $H(\lambda^{(2)}) = \{10,5\}$. We construct a Type II parent $\nu = (4,4,3,1,1,1)$ by first finding the 4-shift of $H(\mu)$, $H(\mu)^{+4} = \{9,6,3,2,1,0\}$, and replacing 0 by 8 to obtain $H(\nu) = \{9,8,6,3,2,1\}$.  
We list out the Type I 8-parents of $\mu = (4,2)$:
\begin{itemize}
    \item $\la^{(1)} = (12,2)$ with $H(\la^{(1)}) = \{13,2\}$ arising via the addition of 8 to the element 5. 
    \item  $\la^{(2)} = (9,5)$ with $H(\la^{(2)}) = \{10,5\}$ arising via the addition of 8 to the element 2.
\end{itemize}
The corresponding Ferrers diagrams are shown below:
\begin{center}
$\sh(\la^{(1)}) = \y{4,2}*[*(lightgray)]{4+8}$ \qquad  $\sh(\la^{(2)}) = \ydiagram{4,2}*[*(lightgray)]{4+5,2+3}$. 
\end{center}
We list out the Type II 8-parents of $\mu = (4,2)$:
\begin{itemize}
    \item $\la^{(3)} = (6,5,3)$ with $H(\la^{(3)}) = \{8,6,3\}$ arising from a 1-shift of $H(\mu)$, and replacing 0 by 8 in $H(\mu)^{+1}$.
    \item $\la^{(4)} = (5,5,3,1)$ with $H(\la^{(4)}) = \{8,7,4,1\}$ arising from a 2-shift of $H(\mu)$, and replacing 0 by 8 in $H(\mu)^{+2}$.
    \item $\la^{(5)} = (4,4,3,1,1,1)$ with $H(\la^{(5)}) = \{9,8,6,3,2,1\}$ arising from a 4-shift of $H(\mu)$, and replacing 0 by 8 in $H(\mu)^{+4}$.
    \end{itemize}
    The corresponding Ferrers diagrams are shown below:
\begin{center}
$\sh(\la^{(3)}) = \ydiagram{4,2}*[*(lightgray)]{4+2,2+3,3}$ \qquad  $\sh(\la^{(4)}) = \ydiagram{4,2}*[*(lightgray)]{4+1,2+3,3,1}$ \qquad  $\sh(\la^{(4)}) = \ydiagram{4,2}*[*(lightgray)]{0,2+2,3,1,1,1}$. 
\end{center}
    \begin{itemize}
    \item $\la^{(6)} = (4,3,3,1,1,1,1)$ with $H(\la^{(6)}) = \{10,8,7,4,3,2,1\}$ arising from a 5-shift of $H(\mu)$, and replacing 0 by 8 in $H(\mu)^{+5}$.
    \item $\la^{(7)} = (4,2,2,1,1,1,1,1,1)$ with $H(\la^{(7)}) = \{12,9,8,6,5,4,3,2,1\}$ arising from a 7-shift of $H(\mu)$, and replacing 0 by 8 in $H(\mu)^{+7}$.
    \item $\la^{(8)} = (4,2,1,1,1,1,1,1,1,1)$ with $H(\la^{(8)}) = \{13,10,8,7,6,5,4,3,2,1\}$ arising from a 8-shift of $H(\mu)$, and replacing 0 by 8 in $H(\mu)^{+8}$.
\end{itemize}
 The corresponding Ferrers diagrams are shown below:
\begin{center}
$\sh(\la^{(6)}) = \ydiagram{4,2}*[*(lightgray)]{0,2+1,3,1,1,1,1}$ \qquad  $\sh(\la^{(7)}) = \ydiagram{4,2}*[*(lightgray)]{0,0,2,1,1,1,1,1,1}$ \qquad  $\sh(\la^{(8)}) = \ydiagram{4,2}*[*(lightgray)]{0,0,1,1,1,1,1,1,1,1}$.
\end{center}
We note that there are no 8-parents arising from the shifts $r_5 = 8 - 5 = 3$ and $r_2 = 8 -2 = 6$ as both $H(\mu)^{+3} = \{8,5,2,1,0\}$ and $H(\mu)^{+6} = \{11,8,5,4,3,2,1,0\}$ contain both 8 and 0 simultaneously.
\end{example}
\cref{parents} can be used to deduce the counts of Type I and Type II $2^R$-parents. Note that there are exactly $\ell(\mu)$ Type I $2^R$-parents as they all arise by choosing a certain element of $H(\mu)$ to be replaced. On the other hand, the values of $r$ for which $H(\mu)^{+r}$ contains $2^R$ (equivalently, have \textbf{no} corresponding $2^R$-parents) are in correspondence with the entries of $H(\mu)$ by letting $r^h = 2^R-h$ for any $h\in H(\mu)$. Thus, there are exactly $2^R-\ell(\mu)$ values of $r$ in the range $1\leq r\leq 2^R$ that lead to a $2^R$-parent of $\mu$. So, by adding the two contributions, we see that any odd-dimensional partition $\mu$ has exactly $2^R$ $2^R$-parents. Iterating this method recovers the count $a(n)$ of odd-dimensional partitions due to McKay \cite{mckay} and Macdonald \cite{mcd}.
\begin{proposition}\label{count}
Let $n = 2^{k_1} + \ldots + 2^{k_r}$ with $k_1 > \ldots > k_r$. Then the number of odd-dimensional partitions of $n$ is given by $a(n) := 2^{k_1 + \ldots + k_r}$.
\end{proposition}
\begin{example}
    If $n = 74 = 2^6 + 2^3 + 2^1$, then the unique odd-dimensional partition of 0 has 2 odd-dimensional 2-parents (\cref{odd-partitions-are-hooks}). Any odd-dimensional partition of 2 has 8 $2^3$-parents giving us a total of 16 odd-dimensional partitions. Then any odd-dimensional partition of 16 has 64 $2^6$-parents, giving us a total of $2^{6+3+1} = 1024$ odd-dimensional partitions of 74.
\end{example}
\subsection{\textnormal{Od} function}\label{odsubs}
To compute $\delta(n) = a_1(n) - a_3(n)$, we subtract the number of partitions with dimension 3 mod 4 from the number of partitions with dimension 1 mod 4. With each partition with dimension 1 modulo 4, we associate the value 1, and with each partition with dimension 3 modulo 4, we associate the value $-1$. The sum of these values over all odd-dimensional partitions of $n$ is $\delta(n)$. We extract these values by using the function $\Od$. 
\begin{definition}
Let $n$ be a positive integer and let $k$ be the largest odd number dividing $n$. We define $\Od(n) = 1$ if $k\equiv 1$ mod 4, and $\Od(n) = -1$ if $k\equiv -1$ mod 4.
\end{definition}
\begin{example}
    Let $n = 1088$ which factors as $64\times 17$. As 17 is the largest number dividing $n$ and $17\equiv 1$ mod 4, we have $\Od(1088) = 1$.  
\end{example}
The condition of $f^\la$ being odd is equivalent to the cancellation of all even factors among the three parts (\cref{rem:three-parts}) of the Frobenius hook-length formula. The next proposition ensures that $\Od(f^\la)$ can be computed as a product of $\Od$ of each of its factors. This method works as $\{-1,1\}$ forms a multiplicative group in $\mathbb{Z}_4$.
\begin{lemma}\label{lem:odmulti}
For all $m, n\geq 1$, $\Od(mn) = \Od(m)\Od(n)$. When $n$ divides $k$, $\Od(k/n) = \Od(k)/\Od(n)$.
\end{lemma}
\begin{proof}
The largest odd number dividing $mn$ is the product of the largest odd numbers dividing $m$ and $n$ respectively. Considering this equality modulo 4 gives us the result. The latter result follows by letting $k = mn$.
\end{proof}
\begin{example}
    Let $m= 36$ and $n = 42$. Then the odd part of $m$ is 9 while the odd part of $n$ is 7. The odd part of $mn = 1512$ is 63. So, we have $\Od(9)\Od(7) = 1\cdot -1 = -1 = \Od(63)$.
\end{example}

\begin{remark}\label{rem:odbinary}
Let $n$ have the binary expansion $[b_kb_{k-1}\ldots b_0]_2$. Let $j$ be the minimum index such that $b_j = 1$. Then, $\Od(n) \equiv [b_{j+1}b_{j}]_2 \mod 4$. Equivalently, one may view the process of finding $\Od(n)$ as starting with the binary expansion of $n$, removing all the rightmost zeros, and returning as the output the rightmost two digits of the newly obtained string mod 4.
\end{remark}

\section{The \workhorse}\label{sec:workhorse}
In this section, we derive a relationship between the values of $\Od(f^\la)$ and $\Od(f^\mu)$ where $\la$ is the $2^R$-parent of $\mu$. As $\Od$ is multiplicative (\cref{lem:odmulti}) and $\Od(n)$ is its own reciprocal for every $n$, it is sufficient for us to consider how the terms in the Frobenius hook-length formula (\cref{HLF}) change in passage from $f^\la$ to $f^\mu$. As two of the three parts (\cref{rem:three-parts}) of the Frobenius hook-length formula contain factorials, we first analyze how factorials behave under $\Od$.

\subsection{Odd part of factorial modulo 4} We derive an explicit formula for $\Od(n!)$ for all non-negative integers $n$. This formula depends on two statistics associated with binary expansions. The first one is the sum of bits $\nu(n)$ in the binary expansion of $n$, or equivalently the number of 1s in the binary expansion of $n$. The other is a statistic $D(n)$ defined in \cite{Bin4} which we now recall.
\begin{definition}
For a non-negative integer $n$, define $D(n)$ to be the number of pairs of consecutive 1s in the binary expansion of $n$. If $n=[b_k\ldots b_0]_2$, then $D(n) = |\{0\leq i \leq k-1: b_i\cdot b_{i+1} =1\}|$.
\end{definition}
\begin{example}
We have $7 = [111]_2$, $10 = [1010]_2$ and $112 = [11011110]_2$.
We compute $D(7) = 2$, $D(10) =0$ and $D(112) =  4$.
\end{example}
We first derive a recursion for the odd part of $n!$ modulo 4. The idea of the proof is to separately analyze the product of odd numbers and even numbers less than $n$, and build the recursion using the fact that $\Od(2k) = \Od(k)$ for all $k\geq 1$.
\begin{lemma}\label{lem:odnfac-rec}
    For any natural number $n$, \[
\Od(n!) = (-1)^{\lfloor (n+1)/4\rfloor}\Od(\lfloor n/2 \rfloor!).
\]
where $\lfloor n \rfloor$ denotes the integral part of $n$.
\end{lemma}
\begin{proof}
    By the multiplicativity of $\Od$ (\cref{lem:odmulti}),
\[
\Od(n!) = \prod\limits_{r=1}^n \Od(r) = \prod\limits_{\substack{1 \leq r \leq n\\ r \text{ odd}}} \Od(r) \prod\limits_{\substack{1 \leq r \leq n\\ r \text{ even}}} \Od(r).
\]
As multiplication by 2 does not affect the odd part, we have $\Od(2r) = \Od(r)$ which gives \[\prod\limits_{\substack{1 \leq r \leq n\\ r \text{ even}}} \Od(r) = \prod\limits_{k=1}^{\lfloor n/2 \rfloor}\Od(2k) =  \prod\limits_{k=1}^{\lfloor n/2 \rfloor}\Od(k) =\Od(\lfloor n/2 \rfloor!).\]
We can split the product over odd numbers depending on their residue modulo 4 as
\begin{align*}
\prod\limits_{\substack{1 \leq r \leq n\\ r \text{ odd}}} \Od(r) &= \prod\limits_{\substack{1 \leq r \leq n\\ r \in 4\mathbb{N}+1}} \Od(r) \prod\limits_{\substack{1 \leq r \leq n\\ r \in 4\mathbb{N}-1}} \Od(r)\\
&= 1\cdot \prod\limits_{\substack{1 \leq r \leq n\\ r \in 4\mathbb{N}-1}} (-1).
\end{align*}

As there are exactly $\lfloor (n+1)/4\rfloor$ terms of the form $4x-1$ less than or equal to $n$ where $x\in \mathbb{N}$, we get

\[
\Od(n!) = (-1)^{\lfloor (n+1)/4\rfloor}\Od(\lfloor n/2 \rfloor!).
\]
\end{proof}

By iterating the above recursion, we can find an explicit formula for $\Od(n!)$.
\begin{lemma}\label{lem:odnfac}
For any natural number $n$, we have \[\Od(n!) = (-1)^{D(n) + \nu(\lfloor n/4 \rfloor)},\] where $\lfloor n \rfloor$ denotes the integral part of $n$.
\end{lemma}
\begin{proof}

Let $n = [b_k\ldots b_0]_2$. In this proof we will omit the brackets and the 2 in the subscript of binary notation, and simply write $b_kb_{k-1}\ldots b_0$ . Taking the integral part of a natural number after division by 2 has the effect of removing the rightmost bit from the binary expansion of $n$, which means $\lfloor \frac{b_\ell\ldots b_0}{2}\rfloor = b_\ell\ldots b_1$. By \cref{lem:odnfac-rec}, we get 
\[
\Od([b_kb_{k-1}\ldots b_0]_2 !) = (-1)^{\lfloor \frac{b_k\ldots b_{0} + 1}{4}\rfloor} \Od([b_kb_{k-1}\ldots b_1]_2 !) 
\]
By another application of \cref{lem:odnfac-rec} on the right, we get 
\[
\Od([b_kb_{k-1}\ldots b_0]_2 !) = (-1)^{\lfloor \frac{b_k\ldots b_{0} + 1}{4}\rfloor + \lfloor \frac{b_k\ldots b_{1} + 1}{4}\rfloor} \Od([b_kb_{k-1}\ldots b_2]_2 !).
\]
Continuing in this manner, we have
\[
\Od(n!) = (-1)^{\lfloor \frac{b_k\ldots b_{0} + 1}{4}\rfloor + \lfloor \frac{b_k\ldots b_{1} + 1}{4}\rfloor + \lfloor \frac{b_k\ldots b_{2} + 1}{4}\rfloor +  \ldots+  \lfloor \frac{b_k b_{k-1} + 1}{4}\rfloor} \Od(b_k!).
\]

Consider the sum in the exponent. For any natural number $m$, we interpret taking the integral part of $m/4$ as removing the {two} rightmost bits in the binary expansion of $m$. Suppose the rightmost pair of bits of $m$ are $b_jb_{j-1}$, then if $0 \leq [b_j b_{j-1}]_2 < 3$, adding 1 to $m$ will only affect the two rightmost bits, thus $\lfloor (m+1)/4\rfloor = \lfloor m/4\rfloor$. On the other hand, if $[b_jb_{j-1}]_2 = 3$, or equivalently, both $b_j$ and $b_{j-1}$ are 1, we have $\lfloor (m+1)/4\rfloor = \lfloor m/4\rfloor + 1$. We deduce that each consecutive pair of 1s appears in the binary expansion of $n$ contributes a $+1$ to the exponent. Thus,
\begin{align*}
    &\left\lfloor \frac{b_k\ldots b_{0} + 1}{4}\right\rfloor + \left\lfloor \frac{b_k\ldots b_{1} + 1}{4}\right\rfloor + \left\lfloor \frac{b_k\ldots b_{2} + 1}{4}\right\rfloor+  \ldots+  \left\lfloor\frac{b_k b_{k-1} + 1}{4}\right\rfloor\\
    &= \left\lfloor \frac{b_k\ldots b_{0} }{4}\right\rfloor + \left\lfloor \frac{b_k\ldots b_{1}}{4}\right\rfloor +   \ldots+ \left\lfloor \frac{b_kb_{k-1}b_{k-2}}{4}\right\rfloor + \left\lfloor \frac{b_k b_{k-1}}{4}\right\rfloor + D(n)
\end{align*}
 Plugging this in the exponent yields
\[
\Od(n!) = (-1)^{D(n)}(-1)^{\lfloor \frac{b_k\ldots b_{0} }{4}\rfloor + \lfloor \frac{b_k\ldots b_{1}}{4}\rfloor +   \ldots+ \lfloor \frac{b_kb_{k-1}b_{k-2}}{4}\rfloor + \lfloor \frac{b_k b_{k-1}}{4}\rfloor} \Od(b_k!).
\]
Using our interpretation of $\lfloor m/4\rfloor$ as removing the two rightmost bits in the binary expansion of $n$, we may write the exponent as:
\[\Od(n!) = (-1)^{D(n)} \cdot (-1)^{b_kb_{k-1}\ldots b_2 + b_{k}b_{k-1}\ldots b_3 + \ldots + b_{k} + 0} \Od(b_k!).
\]
As $b_k \in \{0,1\}$, $b_k! = 1$, and thus $\Od(b_k!) = 1$. So
\[
\Od(n!) = (-1)^{D(n)}(-1)^{{b_k\ldots b_{2} } +  {b_k\ldots b_{3}} +  \ldots+ b_k + 0}.
\]
The sign $(-1)^m$ only depends on the parity of $m$, or in binary expansion notation, $(-1)^{b_k\ldots b_j} = (-1)^{b_j}$. Using this, we obtain,
\[
\Od(n!) = (-1)^{D(n) + b_2 + b_3 + \ldots + b_k}\\
= (-1)^{D(n) + \nu(\lfloor n/4 \rfloor)}.\qedhere
\]
\end{proof}

\subsection{The \workhorse}
In this section, we study the ratio between $\Od(f^\la)$ and $\Od(f^\mu)$ when $\la$ is a $2^R$-parent of $\mu$. From \cref{hookremprop}, we know the removal of a $2^R$-rim hook from the diagram of $\la$ to obtain $\mu$ corresponds to decrementing some entry $h$ in $H(\la)$ by $2^R$ so as to obtain a $\beta$-set of $\mu$. We call the entry $h$ we decrease the \text{affected hook-length}.
\begin{definition}\label{notn:aff}
For partitions $\mu$ and $\la$ such that $\la$ is the $2^R$-parent of $\mu$, define the \textit{affected hook-length} $h^\la_\mu$ to be the element of $H(\la)$ such that replacing $h^\la_\mu$ by $h^\la_\mu - 2^R$ yields a $\beta$-set of $\mu$. In our rim hook removal notation, $H(\mu) \sim_\beta H(\la)[h^\la_\mu \to h^\la_\mu - 2^R]$.
\end{definition}
We see later that the ratio $\Od(f^\lambda)/\Od(f^\mu)$ for $\mu\vdash m$ and its $2^R$-parent $\lambda$ depends on three main factors: the size of $\la$, the affected hook-length, and a statistic $\eta^\lambda_\mu$ computed (with certain caveats) by looking at the number of elements we ``jump over'' when sending $h^\la_\mu$ to $h^\la_\mu-2^R$. Let $H(\la) = \{h_1, h_2, \ldots, h_k\}$ with the entries ordered as $h_1 > h_2 > \ldots > h_k$. There is a $\beta$-set $B$ of $\mu$ containing all the same entries as $H(\la)$ except $h_j = h^\la_\mu$ which is replaced by $h_j - 2^R$.
To motivate the following definition, we emphasize that the only elements that differ between the products $\prod\limits_{\substack{i < j\\ h_i, h_j \in H(\la)}} (h_i - h_j)$ and $\prod\limits_{\substack{i < j\\ h'_i, h'_j \in B}} (h_i' - h_j')$, are the ones involving $h^\la_\mu$.
\begin{definition}\label{defn:inversions}
Let $\lambda$ be a $2^R$-parent of $\mu$ with $|\mu| < 2^R$. Let $h^\lambda_\mu \in H(\lambda)$ be the affected hook-length (\cref{notn:aff}). We define $\eta^\lambda_\mu$ as follows:
\[
\eta^\la_\mu = \left\{x\in H(\la)\backslash \{h^\la_\mu\} :  \dfrac{\Od(|h^\lambda_\mu- x|)}{\Od(|h^\lambda_\mu - 2^R - x|)} = -1\right\}.
\]
\end{definition}
The above definition records the exponent of $-1$ that appears in the product \[ \prod\limits_{\substack{x \in H(\lambda)\\ x \neq h^\lambda_\mu}} \dfrac{\Od(|h^\lambda_\mu- x|)}{\Od(|h^\lambda_\mu - 2^R - x|)},
 = (-1)^{\eta^\la_\mu}.\] In the next subsection, we compute the value of $\eta^\la_\mu$ explicitly by comparing the odd parts modulo 4 of $|h^\la_\mu - x|$ and $|h^\la_\mu - 2^R - x|$. The absolute values are necessary as the elements of a $\beta$-set are always ordered in descending order, and the difference between $x$ and $h^\la_\mu$ (or $h^\la_\mu - 2^R$) must always be positive. To present the \workhorse, we need the next statistic on binary strings.
\begin{definition}
Define $s_2(n)$ to be the sum of the leftmost two bits in the binary expansion of $n$. If $n$ has the binary expansion $[b_kb_{k-1}\ldots b_0]_2$, then $s_2(n)= b_{k-1}+b_k$.
\end{definition}
We are now ready to present the result that forms the crux of our enumeration endeavor. Let $B$ be the $\beta$-set of $\mu$ obtained from $H(\la)$ by decreasing the affected hook-length by $2^R$. We study the effect of $\Od$ on the three parts (\cref{rem:three-parts}) in the product arising from $f^\la/ f^\mu$. We see that $\Od(n!)/\Od(m!)$ yields $(-1)^{s_2(n)}$. The products $\prod_{h\in H(\la)} h!$ and $\prod_{h'\in B} h'!$ differ in exactly one factor, and their ratio is $h^\la_\mu!/ (h^\la_\mu- 2^R)!$. We see that this ratio yields the term $(-1)^{s_2(h^\la_\mu)}$. Finally, we look at the ratio of $\Od$ acting on the Vandermonde products arising from $H(\la)$ and $B$, which by definition and our discussion above, yields $(-1)^{\eta^\la_\mu}$. 
\begin{remark}
    The following proof of \cref{workhorse} showcases the usefulness of our choice of the Frobenius hook-length formula (\cref{HLF}) over the one by Frame-Robinson-Thrall (\cref{eq:frt}), as only need to know the first column hook-lengths and not \text{all} hook-lengths for the $2^R$-core.
\end{remark}

In the following proposition we consider the case of $n\geq 4$. We do so as the removal of the leftmost bit of $n$ only changes the sum of bits of $\lfloor n/4 \rfloor$ only when $n\geq 4$. Our recursion applies to $n\geq 4$, and the smaller remaining values can be dealt with by hand or considered initial values for the recursion.
\begin{proposition}[The \workhorse]\label{workhorse}
Let $n\geq 4$ be a natural number such that $n = 2^R +m$ with $m<2^R$. Let $\mu\vdash m$ and $\lambda\vdash n$ be a $2^R$-parent of $\mu$. Then
\[
\Od(f^{\lambda}) = (-1)^{s_2(n) + s_2(h^\lambda_\mu) + \eta^\lambda_\mu}\Od(f^{\mu}).
\]
\end{proposition}
\begin{proof}
By applying $\Od$ on the formula for $f^\lambda$ in \cref{HLF}, we get
\[
\Od(f^{\lambda}) = (-1)^{D(n) + \nu(\lfloor \frac{n}{4} \rfloor)}\prod\limits_{h\in H(\lambda)} (-1)^{D(h) + \nu(\lfloor \frac{h}{4} \rfloor)} \prod\limits_{\substack{x > y\\ x,y\in H(\lambda)}} \Od(x-y).
\]
Choose the $\beta$-set $B$ of $\mu$ with $|B| = \ell(\lambda)$. Then
\[
\Od(f^{\mu}) = (-1)^{D(m) + \nu(\lfloor \frac{m}{4} \rfloor)}\prod\limits_{h'\in B} (-1)^{D(h') + \nu(\lfloor \frac{h'}{4} \rfloor)} \prod\limits_{\substack{x'> y'\\ x',y'\in B}} \Od(x'-y').\]
We first look at the ratio $\dfrac{(-1)^{D(n) + \nu(\lfloor \frac{n}{4} \rfloor)}}{(-1)^{D(m) + \nu(\lfloor \frac{m}{4} \rfloor)}} = (-1)^{D(n) + \nu(\lfloor \frac{n}{4} \rfloor) - D(n) - \nu(\lfloor \frac{n}{4} \rfloor)}$.
Let $n$ have the binary expansion $[b_R\ldots b_0]_2$ and let $S$ be the largest index smaller than $R$ such that $b_S = 1$. For $m = n-2^R$, we have $m = [b_{S}\ldots b_0]_2$. If $S = R-1$, we lose one consecutive pair of 1s in the binary expansion from $n$ to $m$, while on the other hand, if $S < R-1$, the total number of consecutive pairs of 1s remains unchanged. It follows that $D(n) - D(m) = b_{R-1}$. For the digit sum, we get that removing a 1 from the left decreases the sum by 1, and as $n \geq 4$, the leftmost entry must be at least the third entry from the right. Thus, $\nu(\lfloor n/4 \rfloor) - \nu(\lfloor m/4 \rfloor) =1$. These two facts give us \[(-1)^{D(n) + \nu(\lfloor \frac{n}{4} \rfloor) - D(m) - \nu(\lfloor \frac{m}{4} \rfloor)} = (-1)^{1 +b_{R-1}} = (-1)^{s_2(n)}.\]\\
We have that $h^\lambda_\mu$ is in $H(\lambda)$ but not in $B$ whereas $h^\lambda_\mu-2^R$ is in $B$ but not in $H(\lambda)$. All the other elements of the sets $H(\lambda)$ and $B$ are same. As $n\geq 4$, we have $h^\la_\mu \geq 4$, and so a similar computation as above yields
\begin{align*}
&\dfrac{\prod\limits_{h\in H(\lambda)} (-1)^{D(h) + \nu(\lfloor \frac{h}{4} \rfloor)}}{\prod\limits_{h'\in B} (-1)^{D(h') + \nu(\lfloor \frac{h'}{4} \rfloor)}}\\ &= (-1)^{D(h^\lambda_\mu) + \nu(\lfloor \frac{h^\lambda_\mu}{4} \rfloor) - {D(h^\lambda_\mu - 2^R) - \nu(\lfloor \frac{h^\lambda_\mu - 2^R}{4}\rfloor)}} \\&= (-1)^{s_2(h^\lambda_\mu)}.
\end{align*} 
Finally, we are left with the ratio of products 
\[
\dfrac{\prod\limits_{\substack{x> y\\ x,y\in B}} \Od(x-y)}{\prod\limits_{\substack{x'> y'\\ x',y'\in B}} \Od(x'-y')}
\]
where all the differences not involving $h^\la_\mu$ cancel out. This leaves us with the product \[ \prod\limits_{\substack{x \in H(\lambda)\\ x \neq h^\lambda_\mu}} \dfrac{\Od(|h^\lambda_\mu- x|)}{\Od(|h^\lambda_\mu - 2^R - x|)}.
\] which is exactly $(-1)^{\eta^\la_\mu}$.
Combining all these three results gives us the desired recursion.
\end{proof}
\subsection{Understanding \texorpdfstring{$\eta^\lambda_\mu$}{eta-lambda-mu}} 
Although the recursion in \cref{workhorse} is quite compact, it must be untangled to aid in enumeration. The $s_2$ terms can be computed directly from the binary expansions of $n$ and $h^\la_\mu$, so we must investigate the $\eta^\lambda_\mu$ term. When an entry $h\in H(\lambda)$ is replaced by $h-2^R$, we visualize this as $h$ ``jumping over'' the entries between $h$ and $h-2^R$ where each such jump \textit{may} switch the sign of $\Od$. Some jumps result in a sign change between $\Od(h^\la_\mu)$ and $\Od(h^\la_\mu-2^R)$ , as in the case where $6\in H(\la)$, $h^\la_\mu = 9$ and $H(\mu) = H(\la)[9\to 1]$. Here, $\Od(9-6) = -\Od(6-1)$. However, not all jumps result in a sign change.
 For instance, let $4,8\in H(\lambda)$ and $H(\mu) \sim_\beta H(\lambda)[8\to 0]$. In this case, $\Od(8-4) = \Od(4-0)$ which means that jumping from 8 to 0 over 4 does not change the sign of $\Od$. This phenomenon where some jumps switch the sign while some do not makes analyzing $\eta^\lambda_\mu$ difficult and makes the general problem of computing $\delta(n)$ hard.
\begin{definition}[Indicator function]
Denote by $\mathbb{I}_{X}(x)$, the \textit{indicator function on $X$}, which is equal to 1 if $x\in X$ and equal to 0 otherwise.
\end{definition} 
The following statistic records the number of entries of $H(\la)$ that occur between $h^\la_\mu$ and $h^\la_\mu - 2^R$.
\begin{notation}\label{notation:N-la-mu}
For a partition $\lambda$ of $n = 2^R + m$ with $m < 2^R$ and $h\in H(\lambda)$, let $N_{\lambda}(h)= |\{y \in H(\lambda): h - 2^R< y < h\}|$.
\end{notation}
\begin{remark}\label{rem:height-of-a-hook}
    In terms of rim hooks on diagrams, $N_\la(h^\la_\mu)$ is one less than the number of rows occupied by the rim hook $\la/\mu$. In the literature, this is often called the \textit{height} of the rim hook. For instance, if $H(\la) = \{11, 10, 7, 6, 3, 1\}$ and we obtain $H(\mu) = \{11,7,6,3,2,1\}$ by replacing 10 by 2, then the corresponding $8$-rim hook (in gray) is given by
    \[
    \la/\mu = (6,6,4,4,2,1)/(6,3,3,1,1,1) = \y{6,3,3,1,1,1}*[*(lightgray)]{0,3+3,3+1,1+3,1+1,0}
    \]
    We see that 10 jumps over 7, 6, and 3, which makes $N_\la(10) = 3$, which is exactly one less than the number of rows occupied by the hook.
\end{remark}
\begin{proposition}\label{prop:inversions}
Let $\lambda$ be a $2^R$-parent of $\mu$ with $|\mu| < 2^R$.  Then
\[
\eta^{\lambda}_{\mu} = N_{\lambda}(h^\lambda_\mu) - \mathbb{I}_{H(\lambda)}(h^\lambda_\mu - 2^{R-1}) +  \mathbb{I}_{H(\lambda)}(h^\lambda_\mu + 2^{R-1}) + \mathbb{I}_{H(\lambda)}(h^\lambda_\mu - 2^R-2^{R-1}). 
\]
\end{proposition}
\begin{remark}
   This proposition tells us that to compute $\eta^\lambda_\mu$, we first count the number of entries in $H(\lambda)$ which lie between $h^\lambda_\mu$ and $h^\lambda_\mu - 2^R$. Then, we add or subtract from this quantity depending on if $H(\lambda)$ contains $h^\lambda_\mu + 2^{R-1}$, $h^\lambda_\mu - 2^{R-1}$,  or $h^\lambda_\mu - 2^R - 2^{R-1}$. On a number line, if we consider the interval $[h^\la_\mu-2^R, h^\la_\mu]$, then $h^\la_\mu - 2^{R-1}$ is the midpoint of this interval, and the other two quantities are $2^R$ away on either side of this midpoint. 
\end{remark}
\begin{proof}
In this proof, we denote $h^\lambda_\mu$ simply by $h$, and $x$ is an element of $H(\lambda)$ not equal to $h$. Define the ratio 
\[\rho_h(x) = \dfrac{\Od(|h-x|)}{\Od(|h - 2^R - x|)}.
\] We wish to find all values of $x$ such that $\rho_h(x) = -1$, which allows us to compute $\eta^\la_\mu = \{x \in H(\la): \rho_{h}(x) = -1\}$ by \cref{defn:inversions}. 
%\[(-1)^{\eta^\lambda_\mu} = \prod\limits_{\substack{x\in H(\lambda)\\ x\neq h}} \rho_h(x).\] 
We consider three cases: one where $x$ is smaller than $h-2^R$ and $h$, the other where $x$ is larger than $h$ and $h-2^{R}$, and one where $x$ lies between $h$ and $h-2^R$. In each of these cases, we deal with the absolute value in a different way. 
\begin{enumerate}
\item Let $x < h-2^R < h$. The ratio $\rho_h(x)$ simplifies to $\frac{\Od(h-x)}{\Od(h - 2^R - x)}$. Let $h-x = 2^p\omega$ where $\omega$ is odd and $p$ is a positive integer strictly smaller than $R$. For the numerator of $\rho_h(x)$, we have $\Od(h-x) \equiv \omega \text{ mod } 4$. For the denominator, we find $h-2^R-x = 2^p\omega - 2^R = 2^{p}(\omega - 2^{R-p})$. The part inside the parenthesis is odd as $2^{R-p}$ is even owing to $p<R$. This gives us $\Od(h-x-2^R) \equiv \omega - 2^{R-p}$ mod 4 and we get $\rho_h(x) \equiv \omega \cdot (\omega - 2^{R-p})$ mod 4. If $p<R-1$, $2^{R-p}$ is divisible by 4 which makes $\rho_h(x) \equiv \omega^2$ mod 4, and thus $\rho_h(x) = 1$ as the only odd square residue mod 4 is 1. For $p = R-1$, $\rho_h(x) \equiv \omega( \omega - 2)$ mod 4, which makes $\rho_h(x) =-1$. When $\rho_h(x) = -1$, $h-x = 2^{R-1}\omega$. On one hand, as $x < h-2^R$, we have $h-x > 2^R$ which shows $\omega > 2$. On the other hand $h$ is bounded by $n$ (\cref{max-of-H-lambda}) which is strictly smaller than $2^{R+1}$ giving us $\omega < 4$. The only possible value of $\omega$ is 3. So, $\rho_h(h- 3\cdot 2^{R-1}) = -1$, giving us the term $\mathbb{I}_{H(\lambda)}(h^\lambda_\mu - 2^R - 2^{R-1})$
\item Let $h-2^R < h < x$. This gives $\rho_h(x) = \frac{\Od(x-h)}{\Od(x-h+2^R)}$. Assuming $x - h = 2^p\omega$, and performing a similar analysis as part 1, we obtain $\rho_h(x) \equiv \omega(2^{R-p}+\omega)$ mod 4. For $p < R-1$, we again get $\rho_h(x) = 1$, while for $p = R-1$, we obtain $\rho_h(x) = -1$. The only value of $\omega$ possible in this case is $1$ as $x$ must be smaller than $2^{R-1}$ (\cref{max-of-H-lambda}). So, we get $
\rho_h(h + 2^{R-1}) = -1$, and this gives us the term $\mathbb{I}_{H(\lambda)}(h^\lambda_\mu + 2^{R-1})$.
\item Let $h -2^R <  x < h$. Note that these are exactly the entries enumerated by $N_\la(h)$.
Now, $\rho_h(x)$ becomes $\frac{\Od(h-x)}{\Od(- h +2^R + x)}$. We let $h-x = 2^p\omega$ as before. In this case, we see that the denominator is $\Od(2^R - (h-x)) \equiv 2^{R-p} - \omega$ mod 4 while the numerator is simply $\omega$. This gives us $\rho_h(x) \equiv \omega \cdot (2^{R-p} - \omega)$ mod 4. If $p < R-1$, then $2^{R-p}$ is divisible by 4, and $\rho_h(x) \equiv -\omega^2$ mod 4. Thus all values of $x$ where $p$ is not $R-1$ result in $\rho_h(x) = -1$. If $p = R-1$, we get $\rho_h(x) \equiv \omega(2-\omega)$ mod 4 which is always 1. Thus, for exactly ${N_{\lambda}(h) - \mathbb{I}_{H(\lambda)}(h- 2^{R-1})}$ values, we have $\rho_h(x) = -1$.
\end{enumerate}
\end{proof}

\begin{example}
Pick $\mu = (2,1,1)$ with $H(\mu) = \{4,2,1\}$ and choose the Type II $8$-parent $\lambda = (2,2,2,2,1,1,1,1)$ with $H(\lambda) = \{9,8,7,6,4,3,2,1\}$ which makes $h_\mu^\lambda = 8$. We first compute $N_{\lambda}(h_\mu^\lambda)$ which is $|\{y\in H(\lambda): 0 < y < 8\}| = 6$. Also, $\mathbb{I}_{H(\lambda)}(8 - 4) = 1$ while other indicator functions $\mathbb{I}_{H(\lambda)}(8 + 4)$ and $\mathbb{I}_{H(\lambda)}(8 - 12)$ equal 0. Thus, $\eta_\mu^\lambda = 6 -1  = 5$ in this case. Now consider the Type I $16$-parent $\nu = (10,3,3,3,3,2,2,2)$ of $\lambda$ such that $H(\nu) = \{17,9,8,7,6,4,3,2\}$ and thus $h_\lambda^\nu = 17$. Then, $N_\lambda(h_\lambda^\nu) = |\{y\in H(\nu): 1 < y < 17\}| = 7$. As only $\mathbb{I}_{H(\nu)}(17  - 8) = 1$ we get $\eta_\lambda^\nu = 7 -1 = 6$.
\end{example}

\begin{remark}\label{rem:hook-interpretation}
    We now interpret \cref{prop:inversions} in terms of the rim hook $\gamma = \la/\mu$ on diagrams. These diagrammatic interpretations will not aid us directly in computations, but offer a visual complement to our methods. Recall from \cref{rem:height-of-a-hook} that $N_\la(h^\la_\mu)$ is one less than the number of rows occupied by $\gamma$. Let us denote $h^\la_\mu$ by $h$, and as $h$ is the affected hook-length, $h-2^R\not\in H(\la)$.
    \begin{itemize}
        \item The entry $h + 2^{R-1}$ can only occur in $H(\la)$ in the Type II case as the largest entry in the Type I case is $h$ itself. This means $h+2^{R-1}$ , so $h = 2^{R-1}$. If $h+2^{R-1}$ exists in $H(\la)$ but $h-2^R$ does not, then we may remove a rim hook of size $2^R + 2^{R-1}$ by sending $h+2^{R-1}$ to $h-2^R$. Equivalently, this removal can be interpreted as constructing $\la$ via a two-step process: first, we add a $2^{R-1}$-rim rim hook $\delta$ to some partition $\nu$ to obtain $\mu$, and then add the $2^R$-rim hook $\gamma$ to $\mu$ to obtain $\la$. This is done in such a manner that  $\gamma \cup \delta$ is a rim hook of size $2^R + 2^{R-1}$ added to $\nu$. We can only extend $\gamma$ to $\gamma\cup\delta$ if the southwesternmost cell of $\delta$ lies immediately above the northeasternmost cell of $\gamma$.
        \begin{example}
            In the following example, $\nu = (1,1,1)$, $\mu = (4,2,1)$ is a 4-parent of $\nu$, and $\la = (4,2,2,2,1,1,1,1,1)$ is an 8-parent of $\mu$. We find $H(\nu) = \{3,2,1\}$, $H(\mu) = \{6,3,1\}$, and $H(\la) = \{12,9,8,7,5,4,3,2,1\}$. The following diagram depicts $\nu$ in white, $\delta :=\mu/\nu$ in light gray and $\gamma:= \la/\mu$ in dark gray. Notice how $\gamma$ can be extended to a rim hook of size 12 by combining it with $\delta$.
            \[
            \y{1,1,1}*[*(lightgray)]{1+3,1+1}*[*(gray)]{0,0,1+1,2,1,1,1,1,1}            \]

            As a non-example, we keep the same $\nu$ and $\mu$, but consider $\hat{\la} = (4,3,3,2,1,1,1)$ and notice that $H(\hat{\la}) = \{10,8,7,5,3,2,1\}$ does not contain 12, while also observing that the dark gray $8$-rim hook cannot be extended to a rim hook of size 12.
            \[
            \y{1,1,1}*[*(lightgray)]{1+3,1+1}*[*(gray)]{0,2+1,1+2,2,1,1,1}\]
        \end{example}
        In summary, if our rim hook $\gamma$ can be extended to a rim hook of size $2^R + 2^{R-1}$, then $\mathbb{I}_{H(\la)}(h+2^{R-1}) = 1$.
        
        %%%
        \item The entry $h^\la_\mu - 2^{R-1}$ can only occur in $H(\la)$ when $2^R \leq h^\la_\mu < 2^{R} + 2^{R-1}$. One can prove this condition by assuming the negation that if $h \geq 2^R + 2^{R-1}$ exists in $H(\la)$, then $h - 2^{R-1}\geq 2^R$ exists in $H(\mu)$ where $|\mu| < 2^R$, and thus \cref{max-of-H-lambda} is contradicted. If $h-2^{R-1}$ is not present in $H(\la)$, then we may first send $h$ to $h-2^{R-1}$, and in another step, send $h-2^{R-1}$ to $h-2^R$. This corresponds to decomposing the rim hook $\gamma$ as a disjoint union of two $\gamma'$ and $\gamma''$ where we first remove $\gamma'$ and then $\gamma''$.
        \begin{example}
            Let $\mu = (2,2,1,1)$ with $H(\mu) = \{5,4,2,1\}$. When $\la = (3,3,3,2,2,1)$, we have $H(\la) = \{8,7,6,4,3,1\}$ and the aforementioned decomposition of $\gamma = \la/\mu$ (gray) is not possible as $4\in H(\la)$.
            \[
            \y{2,2,1,1}*[*(gray)]{2+1,2+1,1+2,1+1,2,1}
            \]
            
            On the other hand, when we have the 8-parent $\hat{\la} = (4,3,3,2,2)$ with $H(\hat{\la}) = \{8,6,5,3,2\}$, we can indeed decompose $\hat{\gamma} = \hat{\la}/\mu$ (gray) into $\hat{\gamma}'$ (light gray, top) and $\hat{\gamma}''$ (gray, bottom) as shown
            \[
            \y{2,2,1,1}*[*(gray)]{2+2,2+1,1+2,1+1,2} \to \y{2,2,1,1}*[*(lightgray)]{2+2,2+1,2+1}*[*(gray)]{0,0,1+1,1+1,2}
            \]
        \end{example}
        So, $\mathbb{I}_{H(\la)}(h-2^{R-1}) = 1$ when $\gamma$ \textbf{cannot} be decomposed into two successively removable $2^{R-1}$-rim hooks.

        %%%
        \item The quantity $h - 2^R  - 2^{R-1}$ can only occur in $H(\la)$ when $h > 2^R + 2^{R-1}$, which means $\la$ must be a Type I $2^R$-parent, and so we assume this case. If $h- 2^{R} - 2^{R-1}$ is not an element of $H(\la)$, then we may send $h$ to $h-2^{R} - 2^{R-1}$, which means $\gamma$ can be extended (downwards) to a rim hook $\gamma' = \gamma \cup \delta$ of size $2^{R} + 2^{R-1}$. As the rim hook $\delta$ of size $2^{R-1}$ is added before $\gamma$, the only way for $\gamma$ to not extend to $\gamma'$ is that $\delta$ lies strictly to the south-west of $\gamma$ in $\la$. One may check that if this was not the case, then the removal of $\gamma$ to yield $\mu$ would contain empty boxes in rows above the cells of $\delta$, and thus $\mu$ would not be a partition shape. We can equivalently state this condition as follows: $\gamma$ can be extended to a connected subset of $2^R + 2^{R-1}$ border cells which do not form a hook.
        \boks{0.22}
        \begin{example}
          Let $\mu = (2,2,1,1)$ with $H(\mu) = \{5,4,2,1\}$ where the diagram of $\mu$ is constructed after the addition of a 4-rim hook as $\y{2}*[*(lightgray)]{0,2,1,1}$.
          Let $2^R = 8$, and we consider all four Type I 8-parents of $\mu$. For $\la^1 = (10,2,1,1)$, we have $H(\la^1) = \{13,4,2,1\}$. In this case, both $h^{\la^1}_\mu = 13$ and $h^{\la^1}_\mu-2^R - 2^{R-1} = 1$ occur in $H(\la)$. As we see in the diagram below, the rim hook $\gamma$ (in gray) cannot be extended to a removable rim hook of size 12. We can extend $\la^1/\mu$ to a connected subset of size 12 as shown in light gray, but that is not a rim hook (as its removal does not preserve partition shape). Here,  $\mathbb{I}_{H(\la^1)}(h^{\la^1}_\mu -2^R-2^{R-1}) = 1$ 
          \[
          \sh(\la^1) = \y{2,2,1,1}*[*(gray)]{2+8} \quad \sh(\la^1) = \y{2,2,1,1}*[*(lightgray)]{1+9,2,1}
          \]
          
          For $\la^2 = (9,3,1,1)$ with $H(\la^2)= \{12,5,2,1\}$, $h^{\la^2}_\mu = 12$ exists in $H(\la)$ but $h^{\la^2}_\mu -12 = 0$ does not. As we see from the diagram, $\gamma$ can be extended downwards to a rim hook of size 12, and so $\mathbb{I}_{H(\la^2)}(h^{\la^2}_\mu -2^R-2^{R-1}) = 0$.

          \[
          \sh(\la^2) = \y{2,2,1,1}*[*(gray)]{2+8,2+1} \quad \sh(\la^2) = \y{2,2,1,1}*[*(lightgray)]{2+8,3,1,1}
          \]
          In the cases of $\la^3 = (7,3,3,1)$ with $H(\la^3) = \{10,5,4,1\}$ and $\la^4 = (6,3,3,2)$ with $H(\la^4) = \{9,5,4,2\}$ we have $h^{\la^3}_\mu = 10$ and $h^{\la^4}_\mu = 9$ respectively but as both these are smaller than $12$, the corresponding hooks cannot be extended to border subsets of size 12.
          \boks{0.22}
          \[
           \sh(\la^3) = \y{2,2,1,1}*[*(gray)]{2+5,2+1,1+2} \quad \sh(\la^4) = \y{2,2,1,1}*[*(gray)]{2+4,2+1,1+2,1+1}
          \]
        \end{example}
        \boks{0.3}
        So, $\mathbb{I}_{H(\la)} (h-2^R - 2^{R-1})= 1$ if $\gamma$ can be extended to a non-rim hook border subset of size $2^R + 2^{R-1}$, or equivalently, the unique rim hook of size $2^{R-1}$ lies strictly south-east of $\gamma$.
    \end{itemize}
\end{remark}

\section{Odd-dimensional partitions of Sparse Numbers}\label{sec:proofmain}
In this section, we prove \cref{mainthm} and \cref{maincor} through heavy use of \cref{workhorse} and \cref{prop:inversions}. 

\subsection{The case of sparse numbers}
Recall that a sparse number is a non-negative integer whose binary expansion has no consecutive ones.
As discussed in \cref{ssec:outline}, the case of \textit{sparse numbers} is especially amenable to our enumeration methods. In this section, we assume that $n = 2^{R} + m$ with $m < 2^{R-1}$, which ensures that the binary expansion of $n$ starts with \texttt{10}. This assumption leads to simplified expressions for the \workhorse\, (\cref{cor:sparse-workshorse}) and $\eta^\lambda_\mu$ (\cref{cor:inversions}). Iteratively applying this assumption (for instance, letting $m = 2^S + l$ with $l < 2^{S-1}$) leads to sparse numbers.
\begin{definition}\label{def:sparse}
A positive integer $n$ is called \textit{sparse} if it contains no consecutive 1s in its binary expansion. Equivalently, if $n = [b_k\ldots b_0]_2$, then the product $b_{i+1}\cdot b_i$ equals $0$ for all $i$ such that $0 \leq i < k$.
\end{definition}
\begin{example}
The number $165 = [{10100101}]_2$ is sparse while $91 = [{1011011}]_2$ is not.
\end{example}
We recall the definitions of $a(n)$, $a_i(n)$ and $\delta(n)$ from \cref{ssec:mainresults}. For a non-negative integer $n$, $a(n)$ is the number of odd-dimensional partitions of $n$, $a_i(n)$ is the number of partitions with dimension congruent to $i$ mod 4, and $\delta(n) = a_1(n) - a_3(n)$. Our proof strategy, as discussed in \cref{ssec:outline}, requires constructing partitions iteratively through the addition of hooks. So, as a base case we look at partitions which are a single hook.
\begin{definition}\label{def:hooks}
A partition $\lambda\vdash n$ is called a \textit{hook partition} if it is of the form $(n-b,1, \ldots, 1)$ with $b$ copies of 1 for $0
\leq b \leq n-1$. Equivalently, $\lambda$ is a hook partition of $n$ if $H(\lambda) = \{n, b,b-1, \ldots, 1\}$.
\end{definition}
The dimension of a hook partition $(n-b,1, \ldots, 1)$ is equal to the binomial coefficient $\binom{n-1}{b}$, which follows from \cref{HLF}. The modulo 4 properties of binomial coefficients are discussed in the paper \cite{Bin4} by Davis and Webb. The following proposition in their paper is of use to us.
\begin{proposition}[Thm. 6 in  \cite{Bin4}]\label{davisandwebb}
Let $n\in \mathbb{N}$ be \textbf{not} sparse, then
\[
\Bigg|\Bigg\{k:\binom{n}{k}\equiv 1 \textnormal{ mod } 4\Bigg\}\Bigg| = \Bigg|\Bigg\{k: \binom{n}{k}\equiv 3\textnormal{ mod } 4\Bigg\}\Bigg|.
\]
\end{proposition}
\begin{lemma}\label{lem:hookdims}
We have $\delta(1) = 1$, $\delta(2) = 2$ and $\delta(2^R) = 0$ for $R \geq 2$.
\end{lemma}
\begin{proof}
The cases $\delta(1)$ and $\delta(2)$ can be computed by hand. Now suppose $R \geq 2$. By \cref{odd-partitions-are-hooks}, all odd-dimensional partitions of $2^R$ are hook partitions. Choose $\lambda\vdash 2^R$ such that $H(\lambda) = \{2^R, b, b-1, \ldots, 1\}$. We compute $f^{\lambda} = \binom{2^R-1}{b}$. As $2^R - 1$ contains only 1s in its binary expansion, it is not a sparse number. We now apply \cref{davisandwebb} to deduce that $a_1(2^R) = a_3(2^R)$, which shows $\delta(2^R) = 0$ for $R \geq 2$. 
\end{proof}

For the following exposition, we assume $n = 2^R + m$ with $m < 2^{R-1} < 2^R$. We proceed with the following strategy: choose an odd-dimensional  $2^R$-core $\mu\vdash m$. Consider all partitions of $n$ which are the $2^R$-parents of $\mu$. These are all odd-dimensional partitions by \cref{mcdprop}. By \cref{workhorse}, we can determine $\delta(n)$ in terms of $\delta(m)$. First, we see that \cref{workhorse} simplifies nicely in the $m < 2^{R-1}$ case.
\begin{corollary}[of \cref{workhorse}]\label{cor:sparse-workshorse}
Let $n = 2^R + m$ with $m < 2^{R-1}$. For odd-dimensional partitions, $\lambda \vdash n$ and $\mu\vdash m$ such that $\lambda$ is a $2^R$-parent of $\mu$, 
\[
\Od(f^{\lambda}) = (-1)^{\eta^\lambda_\mu}\Od(f^\mu).
\]
\end{corollary}
\begin{proof}
Suppose $n = [b_Rb_{R-1}\ldots b_0]_2$. As $m$ is strictly smaller than $2^{R-1}$, it must follow that $b_{R-1} = 0$, and the sum of the first two binary digits $s_2(n) = b_R + b_{R-1}$ must be $1$. By \cref{max-of-H-lambda}, we know that the maximum element of $H(\la)$ is at most  $n$, and by the bound on $m$, we deduce $n = 2^R + m < 2^{R} + 2^{R-1}$. Thus all elements of $H(\lambda)$ are strictly smaller than $2^R + 2^{R-1}$. This gives us the inequality $2^R \leq h^\lambda_\mu < 2^R + 2^{R-1}$ which establishes $s_2(h^\lambda_\mu) = 1$. By plugging the values  $s_2(h^\lambda_\mu) = 1$ and  $s_2(n) = 1$ in $
\Od(f^{\lambda}) = (-1)^{s_2(n) + s_2(h^\lambda_\mu) + \eta^\lambda_\mu}\Od(f^{\mu})$ (\cref{workhorse}), we prove the corollary.
\end{proof}

For $m < 2^{R-1}$, the formula for $\eta^\la_\mu$ in \cref{prop:inversions} is simplified.
\begin{corollary}[of \cref{prop:inversions}]\label{cor:inversions}
Let $n = 2^R + m$ with $m < 2^{R-1}$. For odd-dimensional partitions $\lambda \vdash n$ and $\mu\vdash m$ such that $\lambda$ is a $2^R$-parent of $\mu$,
\[
\eta^{\lambda}_{\mu} = N_{\lambda}(h^\lambda_\mu) - \mathbb{I}_{H(\lambda)}(h^\lambda_\mu - 2^{R-1}).
\]
\end{corollary}
\begin{proof}
As established in the proof of \cref{cor:inversions}, $2^R \leq h^\lambda_\mu < 2^R + 2^{R-1}$. On one hand, as $ h^\lambda_\mu - 2^R - 2^{R-1}$ is negative, it cannot appear as a hook-length. Thus, $\mathbb{I}_{H(\lambda)}(h^\lambda_\mu - 2^R - 2^{R-1}) = 0$. On the other hand, $h^\lambda_\mu$ is at least $2^R$ which means $h_\mu^\la + 2^{R-1}$ is at least  $2^R+ 2^{R-1}$, which is strictly larger than $n$. This shows that $2^R + 2^{R-1}$ cannot occur as a hook-length of $H(\la)$ by \cref{max-of-H-lambda}. This means $\mathbb{I}_{H(\lambda)}(h^\lambda_\mu + 2^{R-1}) = 0$. This renders the third and fourth terms of $\eta^{\lambda}_{\mu} = N_{\lambda}(h^\lambda_\mu) - \mathbb{I}_{H(\lambda)}(h^\lambda_\mu - 2^{R-1}) +  \mathbb{I}_{H(\lambda)}(h^\lambda_\mu + 2^{R-1}) + \mathbb{I}_{H(\lambda)}(h^\lambda_\mu - 2^R- 2^{R-1})$ zero, thereby giving us the corollary. 
\end{proof}
\begin{remark}
    We emphasize that the results of \cref{cor:inversions} and \cref{cor:sparse-workshorse} apply to a subset of numbers strictly larger than sparse numbers. For instance, the above corollaries apply to $46 = [101110]_2$ even if it is not a sparse number. It is only when we apply the condition of the binary expansion starting with \texttt{10} iteratively that sparse numbers emerge. Thus in order to state theorems for a general $n$, we also need to know how the recursions simplify when the numbers do not have the binary expansion prefix \texttt{10}.
\end{remark}
\subsection[Counting 2R parents]{Counting $2^R$-parents}\label{ssec:counting-2R-parents}
In this section, we fix a partition $\mu$ and find how many $2^R$-parents of $\mu$ have the same or different dimension modulo 4 compared to $\mu$. We start with $H(\mu)$ and then create $2^R$-parents through the rim hook addition process described in \cref{parents}. This strategy requires us to consider the case of Type I and Type II parents separately. For Type I parents we only replace elements in $H(\mu)$ whereas for Type II parents, we replace elements in $\beta$-sets arising from the $r$-shifts of $H(\mu)$. We establish some notation before we proceed.
\begin{notation}\label{notn:parents}
When $2^R$ is understood, we denote the set of all $2^R$-parents of $\mu$ by $\p(\mu)$. Denote the set of Type I and Type II $2^R$-parents of $\mu$ by $\p_1(\mu)$ and $\p_2(\mu)$ respectively.
\end{notation}
We have $\p(\mu) = \p_1(\mu) \cup \p_2(\mu)$ and from \cref{parents}, we know $|\p_1(\mu)| = \ell(\mu)$ and $|\p_2(\mu)| = 2^R - \ell(\mu)$.
\begin{notation}[Signed-sum]\label{def:signed-sum}
Let $\mu$ be a partition and $\Lambda$ be some set of partitions. Define the \textit{signed-sum} $S_{\Lambda}(\mu) = \frac{1}{\Od(f^\mu)}\sum_{\lambda\in \Lambda} \Od(f^\lambda)$. 
\end{notation}
The quantity $S_{\Lambda}(\mu)$ keeps track of the number of partitions in $\Lambda$ with the same value of $\Od$ as $f^\mu$ minus the number of partitions in $\Lambda$ with a different value of $\Od$ as $f^\mu$. As will become clearer in the subsequent proofs, accounting for the difference is computationally more useful than enumerating partitions with the same and different $\Od$ values separately. This choice of definition computes $\delta(n)$ for $n = 2^R + m$ where $m < 2^R$ using the identity $\delta(n) = \sum\limits_{\la\vdash n} \Od(f^\la)  =\sum\limits_{\mu\vdash m} \Od(f^\mu) S_{\p(\mu)}(\mu)$. We eventually prove the following theorem.
\begin{theorem}\label{thm:parent-sum}
    Let $n = 2^{R}+m$ with $m < 2^{R-1}$. If $\mu$ is a partition of $m$, then \[S_{\p(\mu)}(\mu) = 2 - 2(-1)^m.\]
\end{theorem}

In order to prove \cref{thm:parent-sum}, we split the computations for Type I and Type II $2^R$-parents. In the case of Type II $2^R$-parents, we further consider two cases. The final proof of the theorem emerges in \cref{pf:thm1} after combining the results from \cref{prop:type1}, and \cref{prop:type2}.
%feb 2 2026
\subsubsection{Type I parents}
Recall from \cref{parents} that a Type I $2^R$-parent of $\mu$ is constructed by replacing some $h\in H(\mu)$ by $h^\lambda_\mu = h + 2^R$. Furthermore, $N_\la(h^\la_\mu)$ counts the number of entries strictly between $h$ and $h^\la_\mu$. As we prove in the next proposition, the ratio $\Od(f^\la)/\Od(f^\mu)$ for a Type I $2^R$-parent is entirely dictated by $N_\la(h^\la_\mu)$, which records the number of entries we jump over. Using this, we can compute $S_{\p_1(\mu)}(\mu) = \sum_{\la\in \p_1(\mu)} \Od(f^\la)/\Od(f^\mu)$.

\begin{proposition}\label{prop:type1}
Let $\mu$ be an odd-dimensional partition of $m < 2^{R-1}$. For any $\lambda\in \p_1(\mu)$, $h_\mu^\lambda -  2^{R-1} \not\in H(\lambda)$ and \[S_{\p_1(\mu)}(\mu) = \begin{cases}
    0 & \text{if } \ell(\mu) \text{ is even}\\
    1 & \text{if } \ell(\mu) \text{ is odd}.\\
\end{cases}\]
\end{proposition}
\begin{proof}
Let $H(\mu) = \{h_1, \ldots, h_k\}$. By \cref{parents}, we know that any Type I $2^R$-parent $\la$ of $\mu$ can be constructed by replacing an element $h\in H(\mu)$ by $h^\la_\mu= h+2^R$. Define $\lambda^{(i)} \in \p_1(\mu)$ such that $H(\lambda^{(i)}) = H(\mu)[h_i\to h_i + 2^R]$ for $i$ satisfying $1\leq i \leq k$.

For any $i$, all elements except $h^{\lambda^{(i)}}_\mu$ in $H(\la^{(i)})$ are the same as the elements of $H(\mu)$. If $h^{\la^{(i)}}_\mu - 2^{R-1} \in H(\la^{(i)})$, then $h^{\la^{(i)}}_\mu - 2^{R-1} \in H(\mu)$. As $h^{\la^{(i)}}_\mu > 2^R$, we have $h^{\lambda}_\mu - 2^{R-1} > 2^{R-1} > m$. By \cref{max-of-H-lambda}, this is not possible, and thus $h^{\la^{(i)}}_\mu - 2^{R-1}$ does not occur in $H(\la)$. This, by definition, implies $\mathbb{I}_{H(\lambda)}(h^\lambda_\mu - 2^{R-1}) = 0$ and the formula in \cref{prop:inversions} simplifies to $\eta^{\la^{(i)}}_\mu = N_{\la^{(i)}}(h^{\la^{(i)}}_\mu)$.
Using this, we compute \[
S_{\p_1(\mu)}(\mu) =  \sum\limits_{\lambda\in \p_1(\mu)} \frac{\Od(f^\lambda)}{\Od(f^\mu)}
= \sum\limits_{\lambda\in \p_1(\mu)}(-1)^{\eta^\lambda_\mu}
= \sum\limits_{\lambda\in \p_1(\mu)}(-1)^{N_{\lambda}(h^\lambda_\mu)}.\]
The second equality is established in \cref{cor:sparse-workshorse} while the third equality follows from the immediately preceding discussion.

As all elements of $H(\mu)$ are strictly smaller than $m < 2^{R-1}$, $h^{\lambda^{(i)}}_\mu$ is the largest element in $H(\lambda^{(i)})$. So the elements of $H(\lambda)$ strictly contained between $h_i$ and $h^{\lambda^{(i)}}_\mu = h_i + 2^R$ are exactly the elements in $H(\mu)$ strictly larger than $h_i$. These are exactly the elements $h_1, h_2, \ldots, h_{i-1}$, and thus we compute
\[
N_{\lambda^{(i)}}(h^{\lambda^{(i)}}_\mu)  
=|\{y\in H(\mu): y > h_i\}|= i - 1.\]
This gives us $
S_{\p_1(\mu)}(\mu) = \sum\limits_{i=1}^k (-1)^{i-1}$, which by considering the parity of $\ell(\mu) = \ell(\lambda)$ gives the result.
\end{proof}
\begin{example}
We start with the partition $\mu = (2,1,1)$ such that $H(\mu) = \{4,2,1\}$ with $\Od(f^\mu) = -1$. Its Type I $8$-parents are the partitions $\lambda^{(1)} = (10,1,1)$, $\lambda^{(2)} = (8,3,1)$ and $\lambda^{(3)} = (7,3,2)$ with $\beta$-sets $\{12,2,1\}$, $\{10,4,1\}$, and $\{9,4,2\}$ respectively. To find $\la^{(1)}$, we replace 4 by 12 in $H(\mu)$, and as there are no elements between 4 and 12, we obtain the sign $(-1)^0 = 1$. For $\la^{(2)}$, we replace the 2 by 10 in $H(\mu)$, which requires jumping over the sole element 4, resulting in the sign $(-1)^1$. Computing $\la^{(3)}$ from $\mu$ requires replacing 1 in $H(\mu)$ with 9 which jumps over two entries, namely 2 and 4, thus giving us the sign $(-1)^2$. We compute $S_{\mathcal{P}_1(\mu)}(\mu) = 1 + (-1) + 1 = 1$, which agrees with the proposition as $\ell(\mu) = 3$ is odd.
\end{example}
\begin{remark}
    We can interpret \cref{prop:type1} in terms of hooks by noticing that the $i-1$ in the exponent corresponds to the height of the rim hook (number of rows occupied minus 1). For each row of the diagram of $\mu$, we have a hook whose southwesternmost cell is in that row. The rim hook with its southwesternmost cell in row $i$ contributes the term $(-1)^{i-1}$.
\end{remark}
\subsubsection{Type II parents}
The calculation for Type II $2^R$-parents is more involved. We eventually prove the following proposition:
\begin{proposition}\label{prop:type2}
Let $\mu$ be an odd-dimensional partition of $m< 2^{R-1}$. Then
\[
S_{\p_2(\mu)}(\mu) = \begin{cases}
2 - 2(-1)^{m}, &\text{if } \ell(\mu) \text{ is even}\\
1 - 2(-1)^{m}, & \text{if } \ell(\mu) \text{ is odd}.\\
\end{cases}
\]
\end{proposition}

Fix $\mu\vdash m$ and $n = 2^R + m$ with $m < 2^{R-1}$. Recall that Type II $2^R$-parents are constructed by choosing a $\beta$-set $H(\mu)^{+r}$ of $\mu$ containing 0, and replacing the 0 by $2^R$. In terms of the affected hook-length (\cref{notn:aff}), if $\lambda\in\p_2(\mu)$, then $h^\lambda_\mu = 2^R$. Thus, \cref{cor:inversions} simplifies to 
\begin{equation}\label{eq:etatypeii}
    \eta^\lambda_\mu = N_{\lambda}(2^R) - \mathbb{I}_{H(\lambda)}(2^{R-1}).
\end{equation} The calculation of this quantity comes with the caveat that not all $r$-shifts of $H(\mu)$ have a corresponding $2^R$-parent. For instance, if we want to consider Type II 8-parents of $\mu = (2,1,1)$ with $H(\mu) = \{4,2,1\}$, we need to find the values of $r$ for which the $r$-shifts do not contain an 8. We see that $H(\mu)^{+2} = \{6,4,3,1,0\}$ does not contain an 8, so we can define the corresponding 8-parent $\{6,4,3,1\}$ of $\mu$. We call this 8-parent $\lambda_{[2]}$ where the subscript refers to the 2-shift of the $\beta$-set. On the other hand, $H(\mu)^{+4} = \{8,6,5,3,2,1,0\}$ contains an 8 and thus does not have a corresponding 8-parent. This motivates the following definition.
\begin{notation}\label{notn:valid-values}
For $r \in \{1,\ldots, 2^R\}$, define $\lambda_{[r]}$ to be the partition for which $H(\la_{[r]})$ is constructed by replacing 0 by $2^R$ in $H(\mu)^{+r}$. If $2^R\in H(\mu)^{+r}$, define $\lambda_{[r]} =\mu$. A value of $r$ is called \textit{valid} if $\lambda_{[r]} \neq \mu$.
\end{notation}

We see from \cref{eq:etatypeii}, that the existence of $2^{R-1}$ in $H(\la)$ affects the value of $\eta^\la_\mu$. We can ensure $2^{R-1}$ is present in a $\beta$-set of $\mu$ by shifting $H(\mu)$ by at least $2^{R-1}+1$, as any shift of $r$ introduces the elements $0, 1, \ldots, r-1$ in the $\beta$-set. We denote the set of such $2^R$-parents by $\p_2^{\downarrow}(\mu)$. On the other hand, $H(\mu)^{+r}$ when $1\leq r\leq2^{R-1}$ will contain a $2^{R-1}$ only when there exists an $h$ in $H(\mu)$ such that $h+r = 2^{R-1}$. For each value of $h\in H(\mu)$, we can find $r_h = 2^{R-1}-h$ such that $H(\mu)^{+r_h}$ contains $2^{R-1}$. So exactly $\ell(\mu)$ values of $r$ in $\{1,2, \ldots, 2^{R-1}\}$ will yield $2^{R-1}\in H(\mu)^{+r}$ for $1\leq r\leq 2^{R-1}$. This set of $2^R$-parents arising from valid $r$ values is denoted by $\p_2^{\uparrow}(\mu)$. We describe the two sets in notation, and only keep the parents arising from the valid values of $r$.
\begin{notation}\label{not:typeii}
Let \[\p_2^\uparrow(\mu) := \{\lambda_{[r]} : 1\leq r \leq 2^{R-1}\}\backslash \{\mu\}\] and \[\p_2^\downarrow(\mu) := \{\lambda_{[r]}: 2^{R-1}+1 \leq r \leq 2^{R}\}\backslash \{\mu\}.\]
\end{notation}

We first focus our attention on the $2^R$-parents arising from $r$-shifts with $1\leq r \leq 2^{R-1}$, that is, the set $\p_2^\uparrow(\mu)$. In the proof of the next proposition, we compute $S_{\p_2^\uparrow(\mu)}(\mu)$. The sum $\sum_{\la\in \p_2^\uparrow(\mu)}(-1)^{\eta^\la_\mu}$ eventually simplifies to an expression dependent on the sum $\sum_{h\in H(\mu)} (-1)^h$ which can be computed by only knowing the parity of the elements of $H(\mu)$. Concretely, each even element contributes a $+1$ while an odd element contributes a $-1$. To capture this idea, we define the following notion.
\begin{definition}[Parity gap]
For a finite $X\subset \mathbb{N}$, define the \textit{parity gap} $\g(X)$ to be the number of even elements of $X$ minus the number of odd elements of $X$, i.e., $\g(x) = \sum\limits_{x\in X} (-1)^x$.
\end{definition}
\begin{example}
If $X = \{13, 12, 8, 5, 3, 1, 0\}$ then $\g(X) = 3 - 4 = -1$.
\end{example}
The enumeration result for Type II $2^R$-parents of $\mu$ arising from $r$-shifts with $1\leq r \leq 2^{R-1}$ takes a concise form dependent on the parity of the elements of $H(\mu)$.
\begin{lemma}\label{lem:toprow}
Let $\mu$ be an odd-dimensional partition of $m < 2^{R-1}$, then
\[S_{\p_2^\uparrow(\mu)}(\mu) = 2(-1)^{\ell(\mu)}\g(H(\mu)).\]
\end{lemma}
\begin{proof}
We first show that the largest element in $H(\la_{[r]})$ must be $2^R$. The largest element of $H(\mu)^{+r}$ is $\max(H(\mu))+r$ which itself is \textit{strictly} smaller than $2^R$ as $\max(H(\mu))< m< 2^{R-1}$ (\cref{max-of-H-lambda}) and $r \leq 2^{R-1}$. Thus for $1\leq r \leq 2^{R-1}$, we have $2^R \not\in H(\mu)^{+r}$, which means that each value of $r$ in $\{1,2,\ldots, 2^{R-1}\}$ is valid. So $|\p_2^\uparrow(\mu)| = 2^{R-1}$.
The first column hook-length set of $\lambda_{[r]}$, $H(\la_{[r]})$, is obtained by replacing a 0 in $H(\mu)^{+r}$ by a $2^R$. This makes $2^R$ the largest element of $H(\lambda_{[r]})$ which means that
$N_{\lambda_{[r]}}(2^R) = |\{y\in H(\lambda_{[r])}: 0 <  y < 2^R\}|$ counts every element of $H(\mu)$ except $2^R$ itself. As an $r$-shift introduces $r$ new elements, including a 0 which is replaced by a $2^R$, we find $N_{\lambda_{[r]}}(2^R) = \ell(\mu) + r - 1$. Using \cref{def:signed-sum}, $S_{\p_2^\uparrow(\mu)}(\mu) = \sum_{\la\in \p_2^\uparrow(\mu)}= \frac{\Od(f^\la)}{\Od(f^\mu)}$ which by an application of \cref{cor:sparse-workshorse} becomes $S_{\p_2^\uparrow(\mu)}(\mu) \sum_{\la\in \p_2^\uparrow(\mu)} (-1)^{\eta^\la_\mu}$. A direct application of \cref{eq:etatypeii} replaces $\eta^\la_\mu$  by $N_\la(2^R) - \mathbb{I}_{H(\la)}(2^R)$ to yield
\begin{align*}
S_{\p_2^\uparrow(\mu)}(\mu) &= \sum_{\la\in \p_2^\uparrow(\mu)} (-1)^{N_\la(2^R) - \mathbb{I}_{H(\la)}(2^{R-1})}\\
&=  \sum\limits_{r=1}^{2^{R-1}} (-1)^{\ell(\mu) + r - 1 - \mathbb{I}_{H(\lambda_{[r]})}(2^{R-1})}\\ 
&=(-1)^{\ell(\mu) - 1}\sum\limits_{r=1}^{2^{R-1}} (-1)^{r +\mathbb{I}_{H(\lambda_{[r]})}(2^{R-1})}.
\end{align*}
The second equality follows from the preceding discussion. In the third equality, we factor out $(-1)^{\ell(\mu)-1}$ and replace $-\mathbb{I}_{H(\lambda_{[r]})}(2^{R-1})$ with $+\mathbb{I}_{H(\lambda_{[r]})}(2^{R-1})$ as it preserves the parity of the exponent.
The factor $\mathbb{I}_{H(\lambda_{[r]})}(2^{R-1})$ takes either a value 0 or 1. The value 1 corresponds to the $2^R$-parents $\la$ in $\p_2^\uparrow(\mu)$ where $H(\la)$ contains $2^{R-1}$. As $1\leq r\leq 2^{R-1}$, this can only happen if there exists an $h$ such that $h +r = 2^{R-1}$. In fact, for every $h\in H(\mu)$, we can choose a unique $r_h = 2^{R-1} - h$ which ensures that $2^{R-1} \in H(\mu)^{+r_h}$. This gives rise to an injective map $H(\mu)\rightarrow\p_2^{\uparrow}(\mu)$ given by $h\mapsto \lambda_{[r_h]}$ where the image contains \textit{all} partitions with $2^{R-1}$ as a first column hook-length.  We now break the sum into two parts depending on whether $2^{R-1}$ belongs to $H(\lambda_{[r]})$, or equivalently whether $\mathbb{I}_{H(\lambda_{[r]})}(2^{R-1}) = 1$. 
Let $X = \{r_h : h\in H(\mu)\}$ and $X^c = \{1, \ldots, 2^{R-1}\}\setminus X$.
\begin{align*}
&S_{\p_2^\uparrow(\mu)}(\mu)\\ &= (-1)^{\ell(\mu)-1}\left(\sum\limits_{r_h \in X} (-1)^{r_h +\mathbb{I}_{H(\lambda_{[r_h]})}(2^{R-1})} +\sum\limits_{r\in X^c} (-1)^{r +\mathbb{I}_{H(\lambda_{[r]})}(2^{R-1})}\right)\\
&= (-1)^{\ell(\mu)-1}\left(\sum\limits_{r_h\in X} (-1)^{r_h + 1} +\sum\limits_{r\in X^c} (-1)^{r}\right)
\end{align*}
We add the terms $\sum\limits_{r\in X} (-1)^{r}$  to the second sum $\sum\limits_{r\in X^c} (-1)^{r}$ to write it as a sum over all values of $r$. To ensure we do not change the sum, we also subtract $\sum\limits_{r_h\in X} (-1)^{r_h}$.
\begin{align*}
S_{\p_2^\uparrow(\mu)}(\mu)&= (-1)^{\ell(\mu)-1}\Bigg(\sum\limits_{r_h\in X} (-1)^{r_h + 1} - \sum\limits_{r_h\in X} (-1)^{r_h} 
\\
  & \qquad\qquad\qquad+\sum\limits_{r\in X} (-1)^{r} 
 +\sum\limits_{r\in X^c} (-1)^{r}\Bigg) \\
&= (-1)^{\ell(\mu)-1}\left(\sum\limits_{r_h\in X} \left((-1)^{r_h + 1} - (-1)^{r_h}\right) +\sum\limits_{r=1}^{2^{R-1}} (-1)^{r}\right)
\end{align*}
As $R> 1$, the sum $\sum\limits_{r=1}^{2^{R-1}} (-1)^{r}$ is zero. We also observe that $-(-1)^{r_h} = (-1)^{r_h+1}$, and thus the terms in the first sum are each counted twice. This yields
\begin{align*}
S_{\p_2^\uparrow(\mu)}(\mu)&= 2(-1)^{\ell(\mu)-1}\left(\sum\limits_{r_h\in X} (-1)^{r_h+1}+0\right)\\
&=  2(-1)^{\ell(\mu)}\left(\sum\limits_{r_h\in X} (-1)^{r_h}\right).
\end{align*}
The second equality is obtained by factoring a $(-1)$ from the terms $(-1)^{r_h+1}$ inside the sum and multiplying it with the scalar $(-1)^{\ell(\mu)-1}$ giving $(-1)^{\ell(\mu)}$. Notice that if we have $h\in H(\mu)$ and correspondingly an $r_h$ with $h+r_h = 2^{R-1}$, then $h$ and $r_h$ have the same parity. Thus we can rewrite the sum over $r_h\in X$ as a sum over $h\in H(\mu)$. By using the definition of the parity gap, $\g(H(\mu))$,
\[
S_{\p_2^\uparrow(\mu)}(\mu)=2(-1)^{\ell(\mu)}\left(\sum\limits_{\substack{h\in H(\mu)}} (-1)^{h}\right)=
2(-1)^{\ell(\mu)}\g(H(\mu)).\qedhere
\]
\end{proof}
One may wonder why the above statement was relegated to being a lemma rather than a full-fledged proposition. Our work here is not done as we can further analyze the outputs for $\g(H(\mu))$. We emphasize the following lemma holds for all odd-dimensional partitions and the sizes are \textbf{not} subject to the constraint $m < 2^{R-1}$.
%Feb 3 2026
\begin{lemma}\label{prop:parity}
For an odd-dimensional partition $\lambda$,
\[
\g(H(\lambda)) = \begin{cases}
0, &\text{if }\ell(\lambda) \text{ is even and } n \text{ is even}\\
1, &\text{if }\ell(\lambda) \text{ is odd and } n \text{ is even}\\
2, &\text{if }\ell(\lambda) \text{ is even and } n \text{ is odd}\\
-1, &\text{if }\ell(\lambda) \text{ is odd and } n \text{ is odd}\\

\end{cases}
\]
or more concisely,
\[
\g(H(\lambda)) = \begin{cases}
1 - (-1)^n, &\text{if }\ell(\lambda) \text{ is even}\\
(-1)^n, & \text{if }\ell(\lambda) \text{ is odd}.\\
\end{cases}
\]
\end{lemma}

\begin{proof}
Let $\la$ be a $2^R$-parent of $\mu$. Then by \cref{parents}, there exists an $s$ satisfying $0\leq s\leq 2^R$ and an $h\in H(\mu)^{+s}$ such that $H(\lambda) \sim_\beta H(\mu)^{+s}[h\to h+2^R]$. In the case of Type I $2^R$-parents, we replace some $h\in H(\mu)$ by $h+2^R$, while in the case of Type II $2^R$-parents, we consider an $s$-shift of $H(\mu)$ and replace 0 by $2^R$. Thus, the only difference between the sets $H(\la)$ and $H(\mu)^{+s}$ is that $h$ in $H(\mu)^{+s}$ is replaced with $h+2^R$ in $H(\la)$. As $h$ and $h+2^R$ have the same parity, replacing $h$ by $h+2^R$ preserves the number of elements of a given parity, and thus $\g(H(\lambda)) = \g(H(\mu)^{+s})$.
As $X^{+1} = \{x+1\mid x\in X\}\cup\{0\}$, we see that under a 1-shift all odd entries become even, all even entries become odd and an additional even entry $0$ is introduced. Thus, $\g(X^{+1}) = 1 - \g(X)$ where $- \g(X)$ arises from swapping the parity of entries and the $+1$ corresponds to the newly introduced zero. Now note that $\g(X^{+2}) = 1-\g(X^{+1})= 1- (1-\g(x)) = \g(x)$. Continuing this, we get
\[
\g(X^{+s}) = 
\begin{cases}
\g(X), &\text{if } s \text{ is even}\\
1-\g(X),&\text{if }s \text{ is odd.}
\end{cases}
\]
Assume $n = 2^{b_l} + \ldots + 2^{b_1} + |\alpha|$ with $b_l > \ldots > b_1>1$ and $\alpha = \varnothing$ or $(1)$.
From \cref{parents} and \cref{mcdprop}, any odd-dimensional partition $\lambda$ of $n$ can be constructed successively as
\[
\alpha =\nu^0 \subset \nu^{1} \subset \ldots \nu^{l}=\lambda
\]
where for each $1\leq i \leq l$, $\nu^i/\nu^{i-1}$ is a $2^{b_i}$-rim hook. In terms of $\beta$-sets, this means that we can find $s_i$ and $h_i\in H(\nu^{i-1})^{+s_i}$ such that $H(\nu^i) = H(\nu^{i-1})^{+s_i}[h_i \to h_i+2^{b_i}]$. For all $i$, the substitution of $h_i$ by $ h_i+2^{b_i}$ does not affect the number of elements in $H(\nu^i)$ and $H(\nu^{i-1})$ with the same parity. We can write 
\[
\g(H(\lambda)) = \g(X^{+\sum\limits_{i=1}^l s_i})
\]
where $\sum\limits_{i=1}^l s_i$ is number of new entries with respect to $H(\alpha)$ and is thus equal to $\ell(\lambda) -\ell(\alpha)$. So, 
\[\g(H(\lambda)) = 
\begin{cases}
\g(H(\alpha)), & \text{if }\ell(\lambda) -\ell(\alpha)\text{ is even}\\
1 -\g(H(\alpha)),&\text{if }\ell(\lambda) -\ell(\alpha)\text{ is odd.}
\end{cases}
\]
When $\alpha = \varnothing$, $\g(H(\alpha)) = 0$ and $\ell(\alpha) = 0$. In this case, $n$ is even and we obtain
\[\g(H(\lambda)) = 
\begin{cases}
0 & \ell(\lambda) \text{ is even}\\
1 &\ell(\lambda)\text{ is odd.}
\end{cases}
\]
When $\alpha = (1)$, we have $\g(\{1\}) = -1$ and $\ell(\alpha) = 1$. This is the case when $n$ is odd, and
\[\g(H(\lambda)) = 
\begin{cases}
2, &\text{if }\ell(\lambda)\text{ is even}\\
-1, & \text{if }\ell(\lambda) \text{ is odd.}
\end{cases}
\]
The concise formula for $\g(H(\la))$ follows from these computations.
\end{proof}
\begin{example}
Let $\lambda$, $\mu$, $\pi$ and $\kappa$ be partitions such that
\begin{itemize}
\item $H(\mu) = \{1\}$ 
\item $H(\la) = H(\mu)^{+2}[0\to 4] =  \{4,3,1\}$
\item $H(\pi) = H(\la)^{+6}[0\to 8] = \{10,9,8,7,5,4,3,2,1\}$
\item $H(\kappa) = H(\pi)[5\to 21] = \{21,10,9,8,7,4,3,2,1\}$.
\end{itemize}
The partition $\kappa$ is a partition of $16 + 8 + 4 + 1 = 29$ which is an odd number. The length of $\kappa$ is 9, i.e. odd, as there are 10 entries in $H(\kappa)$. By \cref{prop:parity}, we know that the number of odd numbers must exceed the number of even numbers by 1 in $H(\kappa)$. A quick glance at the entries of $H(\kappa)$ tells us that this is indeed true.
\end{example}

We can now state the result for $\p_2^\uparrow(\mu)$ in terms of $m$ with dependence on the length of $\mu$.
%Mar 5
\begin{proposition}\label{prop:toprow}
    Let $\mu$ be an odd-dimensional partition of $m < 2^{R-1}$ such that $R\geq 2$. Then
 \[
S_{\p_2^\uparrow(\mu)}(\mu) = \begin{cases}
2(1 + (-1)^{m+1}), &\text{if } \ell(\mu) \text{ is even}\\
2(-1)^{m+1}, & \text{if } \ell(\mu) \text{ is odd}.\\
\end{cases}
\]
\end{proposition}
\begin{proof}
    From \cref{lem:toprow}, we have the formula $S_{\p_2^\uparrow(\mu)}(\mu) = 2(-1)^{\ell(\mu)}\g(H(\mu))$. Now we use the result from \cref{prop:parity}. 
    When $\ell(\mu)$ is even, $(-1)^{\ell(\mu)} = 1$ and $\g(H(\mu))$ is $1 - (-1)^m = 1 + (-1)^{m+1}$. This shows that when $\ell(\mu)$ is even, $\g(H(\mu)) = 2(1 + (-1)^{m+1})$ as stated. For the case when $\ell(\mu)$ is odd, we have $(-1)^{\ell(\mu)} = -1$ and $\g(H(\mu)) = (-1)^m$. This gives us $\g(H(\mu)) = 2(-1)(-1)^m = 2(-1)^{m+1}$ for partitions $\mu$ with odd length.
\end{proof}

We now turn our attention to the values of $r$ between $2^{R-1}+1$ and $2^R$. In this case, we see not all values of $r$ are valid, and thus have no corresponding $2^R$-parents (\cref{notn:valid-values}). In fact, we can show that the values of $r$ that are \textbf{not} valid are exactly $2^{R-1}$ more than the values $r_h = 2^{R-1}-h$ as defined in the proof of \cref{lem:toprow}. As we saw in that proof, these $r_h$-shift values, and by transitivity the invalid $r$ values in $\{2^{R-1}, \ldots, 2^{R}\}$, are in bijection with the elements of $H(\mu)$. The proof of the following proposition is involved but we provide a broad idea through the following hypothetical example. Recall that each $2^R$-parent $\la$ has an associated sign given by $\Od(f^\la)/\Od(f^\mu)$. For $R = 4$, $r$ may inclusively take a value between 9 and 16. Recall from \cref{notn:valid-values} that if $r$ is a valid value, then $\la_{[r]}$ is the $2^R$-parent arising from the $r$-shift by replacing $0$ in $H(\mu)^{+r}$ by $2^R$, and $\la_{[r]} = \mu$ otherwise. Suppose $\la_{[10]}$ has the sign $-1$, and 11, 12, 13 are invalid $r$ values. If 14 is a valid shift value, then we claim that the sign associated with $\la_{[14]}$ is $+1$, that is, the opposite of the immediately preceding $2^R$-parent. We show in general that the signs of $2^R$-parents will alternate for consecutive valid shift values (while skipping the invalid values). So as an example if the valid shift values are $9, 10, 14, 16$, and the sign of $\la_{[9]}$ is $+1$, then the sign of $\la_{[10]}$ is $-1$, the sign of $\la_{[14]}$ is $+1$ and the sign of $\la_{[16]}$  is $-1$. 
\begin{proposition}\label{lem:bottomrow}
Let $\mu$ be a partition of $m < 2^{R-1}$.
For $2^{R-1}+ 1 \leq r \leq 2^R$, we have $\lambda_{[r]} \in \p_2^\downarrow(\mu)$ if and only if $2^{R-1}\not\in H(\lambda_{[r-2^{R-1}]})$. Furthermore, \vspace*{-5pt}
\[
S_{\p_2^\downarrow(\mu)}(\mu) = \begin{cases}
0, &\text{if }\ell(\mu) \text{ is even}\\
1, &\text{if } \ell(\mu) \text{ is odd}.\\
\end{cases}
\]
\end{proposition}
\begin{proof}
\begin{comment}
If $2^{R-1}$ appears in $H(\lambda_{[r-2^{R-1}]})$, then there exists an $h\in H(\mu)$ such that $h + r - 2^{R-1} = 2^{R-1}$. So, $h + r = 2^R$ which shows that $2^R \in H(\mu)^{+r}$.
So, by running the argument in reverse, we get that $2^{R-1}\in H(\lambda_{[r-2^{R-1}]})$ if and only if $2^R \in H(\mu)^{+r}$ for $2^{R-1}+1 \leq r \leq 2^R$.
\end{comment}
For each $h\in H(\mu)$ we can find an $r$ such that $h + r = 2^R$. The values of $r$ obtained in this manner are the only invalid values of $r$. Subtracting $2^{R-1}$ on both sides of $h + r = 2^R$ gives the equality $h+ (r - 2^{R-1}) = 2^{R-1}$ from which the claim $\lambda_{[r]} \in \p_2^\downarrow(\mu)$ if and only if $2^{R-1}\not\in H(\lambda_{[r-2^{R-1}]})$ follows. Now we compute $S_{\p_2^\downarrow(\mu)}(\mu)$.

From \cref{cor:inversions} and \cref{eq:etatypeii}, we have $\eta^\lambda_\mu = N_{\lambda}(2^R) - \mathbb{I}_{H(\lambda)}(2^{R-1})$. By noticing that a shift of $r\geq 2^{R-1}$ would force the $\beta$-set $H(\mu)^{+r}$ to always contain a $2^{R-1}$, $\eta^\la_\mu= N_{\lambda_{[r]}}(2^R) - 1$. We now compute $N_{\lambda_{[r]}}(2^R)$, that is, the number of entries in $H(\la_{[r]})$ strictly less than $2^R$. Recall from \cref{notn:valid-values} that $\la_{[r]}$ is the $2^R$-parent arising by replacing $0$ in $H(\mu)^{+r}$ by $2^R$.
%notations are immediately recalled when used
Suppose $\lambda_{[r]}, \lambda_{[r+s]}\in \p_2^\downarrow(\mu)$ and $\lambda_{[r+1]} = \ldots = \lambda_{[r+s-1]} = \mu\not\in \p_2^\downarrow(\lambda)$. In other words, $r$ and $r+s$ are valid values of $r$-shifts while $r+1, r+2, \ldots, r+s-1$ are not. From our discussion above, a shift $r+i$ for $1\leq i \leq s-1$ can be invalid only if there exists some element $h\in H(\mu)$ such that $h + r+i \in H(\mu)^{r+i}$ is equal to $2^R$. Call the element $h_{k+1}$ that satisfies $h_{k+1} + r+1 = 2^R$. In order for every successive $k$-shift from $k = r+1$ to $k = r+s-1$ to be invalid, the consecutive integers $h_{k+1} -1, h_{k+1} - 2,\ldots, h_{k+1} - (s-2)$ must be elements of $H(\mu)$. For conciseness, we define $h_{k+i} = h_{k+1} - i + 1$ for $1\leq i\leq s-1$. So our first column hook-length set for $\mu$ takes the form \[H(\mu) = \{h_1, \ldots, h_k, h_{k+1}, \ldots, h_{k+s-1},\ldots, h_{l}\}\] where $h_{k+i} + r + i = 2^R$ by definition. 
Note that $h_k + r >2^R$. We cannot have $h_k + r = 2^R$ as $r$ is a valid shift, and if $h_k + r < 2^R$, then $h_{k+1} + r + 1 < h_k + r + 1$ cannot equal $2^R$. Similarly, one may argue that $h_{k+s} + r + s < 2^R$.

Now in order to compute $N_{\la_{[x]}}(2^R)$, the number of entries strictly smaller than $2^R$ in $H(\la_{[x]})$, we need to disregard all entries greater than or equal to $2^R$. We define a quantity to keep track of these larger entries.
Let $j_x$ be the number of elements strictly bigger than $2^R$ in $H(\lambda_{[x]})$.

When we $x$-shift a $\beta$-set, we add exactly $x$ new entries, one of which is zero. If we disregard all entries greater than $2^R$ as well as the zero (which will be replaced by $2^R$ in passage to $\la_{[x]}$), we get
$N_{\lambda_{[x]}}(2^R) = \ell(\mu) + x - 1 - j_x$. If we shift by $r+s$, we see 
$h_{k+i} + r + s > 2^R$ for $1\leq i \leq s-1$, and thus $j_{r+s} = j_r + s - 1$ as exactly $s-1$ values exceed $2^R$ in addition to the values already exceeding $2^R$ in $H(\la_{[r]})$. By plugging $x = r+s$ in the expression for $N_{\lambda_{[x]}}(2^R)$, we get $N_{\lambda_{[r+s]}}(2^R) = \ell(\mu) + (r+s) - 1 - j_{r+s}$ which after using the relationship we just established between $j_{r+s}$ and $j_{r}$, we get $N_{\lambda_{[r+s]}}(2^R) = \ell(\mu) + (r + s - 1) - (j_r + s - 1)$ which simplifies to $\ell(\mu) + r - j_r$. Comparing this to the formula $ N_{\lambda_{[x]}}(2^R)$ for $x = r$, we see that $N_{\lambda_{[r+s]}}(2^R) =N_{\lambda_{[r]}}(2^R) + 1$. This shows that for each consecutive $2^R$-parent of $\mu$ (when the invalid values are ignored), the value of $\eta_\mu^\lambda$ increases by 1.

Let $r_1 < \ldots < r_{2^{R-1} - \ell(\mu)}$ be the valid values of $r$-shifts, that is, $\{\lambda_{[r_1]}, \ldots, \lambda_{[r_{2^{R-1} - \ell(\mu)}]}\}= \p_2^\downarrow(\mu)$. By the above discussion, we deduce $N_{\lambda_{[r_i]}}(2^R) = N_{\lambda_{[r_1]}}(2^R) + i - 1$. Furthermore, we can express $N_{\lambda_{[r_i]}}(2^R)$ in terms of $j_{r_1}$ as \[N_{\lambda_{[r_i]}}(2^R) = (\ell(\mu) + r_1 - 1 - j_{r_1}) + i - 1 = \ell(\mu) + r_1 +i - j_{r_1} - 2.\] Suppose $r_1 = 2^{R-1} + z$ for some $1\leq z\leq 2^{R-1}$, then as $r_1$ is the smallest valid value, any shift value before this must be invalid. This means that $j_{r_1} = z-1$ as each set $H(\mu)^{+(2^{R-1} + 1)}, \ldots,  H(\mu)^{+(2^{R-1} + z - 1)}$ will contain $2^R$. We obtain \[N_{\lambda_{[r_i]}}(2^R) = \ell(\mu) + (2^{R-1}+ z) + i - (z-1) - 2 = \ell(\mu) + 2^{R-1} + i -1.\]
\vspace{0.2cm}

Using \cref{def:signed-sum}, \cref{cor:sparse-workshorse} and \cref{eq:etatypeii}, we compute
\begin{align*}
S_{\p_2^\downarrow(\mu)}(\mu) &= \sum_{\la\in \p_2^\downarrow(\mu)} (-1)^{\eta^\la_\mu}\\  
&=\sum\limits_{i = 1}^{2^{R-1}-\ell(\mu)} (-1)^{N_{\lambda_{[r_i]}}(2^R) - \mathbb{I}_{H(\la)}(2^{R-1})}\\
&=\sum\limits_{i = 1}^{2^{R-1}-\ell(\mu)} (-1)^{N_{\lambda_{[r_i]}}(2^R) - 1}\\
&= \sum\limits_{i = 1}^{2^{R-1}-\ell(\mu)}(-1)^{(\ell(\mu) + 2^{R-1} + i -1)-1}
\end{align*}
In the exponent of $(-1)$, the terms $2^{R-1} - 2$ do not contribute as they are even owing to $R\geq 2$. We factor out $(-1)^{\ell(\mu)}$ to obtain  $(-1)^{\ell(\mu)}\sum\limits_{i = 1}^{2^{R-1}-\ell(\mu)}(-1)^{i} $. Looking at the sum, $\ell(\mu)$ is even, then the terms in the sum pair up and cancel out, while if $\ell(\mu)$ is odd, we get one positive term $+1$. This shows the statement of the proposition.
\end{proof}
\begin{remark}
  \ytableausetup{aligntableaux = top}
    We can understand the claim that $\Od(h^{\la_{[r]}}_\mu)$ and $\Od(h^{\la_{[r+s]}}_\mu)$ \cref{lem:bottomrow} have opposite signs in terms of hooks on diagrams. Let $\mu = (2,1,1,1,1)$ and consider its 8-parents $\la_{[3]}$ and $\la_{[8]}$. We see that \[
    \sh(\la_{[8]}) = \y{2,1,1,1,1}*[*(gray)]{0,0,0,0,0,1,1,1,1,1,1,1,1}\quad \sh(\la_{[3]}) = \y{2,1,1,1,1}*[*(gray)]{0,1+1,1+1,1+1,1+1,2,1,1}.\] The rim hook in $\la_{[8]}$ occupies 8 rows while the rim hook in $\la_{[3]}$ occupies 7. The following diagrams correspond to the invalid shifts 7, 6, 5, 4 respectively from left to right:
    \[
    \y{2,1,1,1,1}*[*(gray)]{0,0,0,0,0,2,1,1,1,1,1,1} \quad \y{2,1,1,1,1}*[*(gray)]{0,0,0,0,1+1,2,1,1,1,1,1} \quad \y{2,1,1,1,1}*[*(gray)]{0,0,0,1+1,1+1,2,1,1,1,1}  \quad \y{2,1,1,1,1}*[*(gray)]{0,0,1+1,1+1,1+1,2,1,1,1}  
    \]
    As the shifts decrease by 1 at each step from 7 to 3, the bottommost cell of the rim hook goes up by one row and so does the topmost cell of the hook, thus maintaining the numbers of rows occupied as 7. Comparing an 8-shift to a 7-shift, however, we see that the bottom-most cell goes up one row, but the top cell remains in the same row, which decreases the number of rows occupied by 1.  
    This is the diagrammatic reason why the sign flips between successive parents.
\end{remark}
The proof of \cref{prop:type2} follows from \cref{prop:toprow} and \cref{lem:bottomrow}.
\begin{proof}[Proof of \cref{prop:type2}]
From \cref{prop:toprow},
\[
S_{\p_2^\uparrow(\mu)}(\mu)
= \begin{cases}
2 - 2(-1)^m, &\text{if }\ell(\mu) \text{ is even}\\
2(-1)^{m+1}, & \text{if }\ell(\mu) \text{ is odd},\\
\end{cases}
\] and 
from \cref{lem:bottomrow},
\[
S_{\p_2^\downarrow(\mu)}(\mu) = \begin{cases}
0, &\text{if }\ell(\mu) \text{ is even}\\
1, &\text{if } \ell(\mu) \text{ is odd}.\\
\end{cases}
\]
As $S_{\p_2(\mu)}(\mu) = S_{\p_2^\uparrow(\mu)}(\mu) + S_{\p_2^\downarrow(\mu)}(\mu)$, we add the result of \cref{lem:toprow} and \cref{lem:bottomrow} to prove the claim in \cref{prop:type2}
\[
S_{\p_2(\mu)}(\mu) = \begin{cases}
2 - 2(-1)^{m}, &\text{if }\ell(\mu) \text{ is even}\\
1 - 2(-1)^{m}, & \text{if }\ell(\mu) \text{ is odd}.\qedhere\\
\end{cases}
\]
\end{proof}

\subsection{Proof of the main theorems}
We first establish \cref{thm:parent-sum} and then use it to prove \cref{mainthm}.
\begin{proof}[Proof of \cref{thm:parent-sum}]\label{pf:thm1}
Recall that $\p(\mu)$ denotes the set of $2^R$-parents of $\mu$ and $S_{\p(\mu)}(\mu) = S_{\p_1(\mu)}(\mu) + S_{\p_2(\mu)}(\mu)$. The result of \cref{prop:type1}
\[S_{\p_1(\mu)}(\mu) = \begin{cases}
    0 & \text{if } \ell(\mu) \text{is even}\\
    1 & \text{if } \ell(\mu) \text{is odd},\\
\end{cases}\]
 when added to the result of \cref{prop:type2}
 \[
S_{\p_2(\mu)}(\mu) = \begin{cases}
2 - 2(-1)^{m}, &\text{if }\ell(\mu) \text{ is even}\\
1 - 2(-1)^{m}, & \text{if }\ell(\mu) \text{ is odd},
\end{cases}
\]
yields $S_{\p(\mu)}(\mu) = 2 - 2(-1)^m$. 
\end{proof}

We are now ready to prove the main theorem of this paper that establishes the recursion for $\delta(n)$ when $n = 2^{R}+m$ with $m < 2^{R-1}$.
\begin{proof}[Proof of \cref{mainthm}]
To compute $\delta(n) = a_1(n) - a_3(n)$, we associate $+1$ to each partition with dimension 1 modulo 4, and $-1$ to each partition with dimension 3 modulo 4. This gives
\[
\delta(n) = \sum\limits_{\substack{\\\lambda\vdash n\\f^\lambda \text{ is odd}}} \Od(f^\lambda)
\]
We have $S_{\p(\mu)} = \sum_{\la\in \p(\mu)} \frac{\Od(f^\la)}{\Od(f^\mu)}$, or equivalently $S_{\p(\mu)} \Od(f^\mu) = \sum_{\la\in \p(\mu)} \Od(f^\la)$ for each odd-dimensional partition $\mu$ of $m$. As each odd-dimensional $\la\vdash 2^R + m$ appears as a parent of a unique odd-dimensional $\mu\vdash m$, we can write
\[
\delta(n) = \sum\limits_{\substack{\mu\vdash m\\f^\mu \text{ is odd}}} S_{\p(\mu)}({\mu})\Od(f^{\mu})
\]
By using the result $S_{\p(\mu)}(\mu) = 2 - 2(-1)^m$, we get 
\[
\delta(n) = (2 - 2(-1)^m) \sum\limits_{\substack{\mu\vdash m\\f^\mu \text{ is odd}}} \Od(f^{\mu}).
\]

The sum above, in analogy with $\delta(n)$, is equal to $\delta(m)$. So,
\[
\delta(n) = (2-2(-1)^m )\delta(m).
\]
By the observation that $n$ and $m$ have the same parity, we find for $n \geq 4$
\[
\delta(n) = \begin{cases}
0, & \text{if }n \text{ is even}\\
4 \delta(m), &\text{if } n \text{ is odd}. \qedhere
\end{cases}
\]
\end{proof}
\begin{proof}[Proof of \cref{maincor}]
Let $n = 2^{k_1} + 2^{k_2} + \ldots + 2^{k_l}$ with $k_1 > k_2 > \ldots > k_l$ be a sparse number; which means $k_i > k_{i+1} + 1$ for all $1\leq i < l$.
Repeatedly applying \cref{mainthm}, we obtain 
\[
\delta(n) = \delta(n-2^{k_1})\cdot \delta(n-2^{k_1}- 2^{k_2}) \cdot \ldots \cdot \delta(2^{k_l})
\]
When $n$ is even, $k_l >0$. By \cref{lem:hookdims}, $\delta(2^{k_l}) =0$ and thus $\delta(n) = 0$. When $n$ is odd, $k_l = 0$. In this case, $\delta(2^{k_l}) = \delta(1) = 1$. As each bit except the rightmost one contributes a factor of 4, we obtain
\[
\delta(n) = \begin{cases}
2\text{,} & \text{if }n = 2\\
0, &\text{if } n > 2 \text{ is even}\\
4^{\nu(n) - 1}, & \text{if } n\text{ is odd}.\qquad\qquad\qedhere
\end{cases}
\]
\end{proof} 

%Feb 4 2026
%Mar 1 2026
%May 23, 2026
\section{Odd-dimensional partitions of \texorpdfstring{$n = 2^R + 2^{R-1}$}{n = 2R + 2R-1}}\label{sec:proof11}
We extend the results of the previous section beyond the sparse case. We consider positive integers $n$ of the form $2^R + 2^{R-1}$ with $R\geq 1$. Equivalently, these are numbers $n$ whose binary expansions are of the form $[{1100...}]_2$, i.e., $\nu(n) = 2$ and $D(n) = 1$. We now prove \cref{11thm} which states that for $n = 2^R + 2^{R-1}$ with $R\geq 1$,
\[
\delta(n) = \begin{cases}
2, & \text{if } R = 1\\
8, & \text{if } R = 2\\
0, & \text{otherwise}.
\end{cases}
\]
We use similar methods as we did in the previous section but with slight modifications. We state a corollary of \cref{workhorse} which is pertinent in this case. Throughout this section, we assume  $R\geq 2$ as \cref{workhorse} requires $n\geq 4$. The case of $R = 1$ is straightforward and can be dealt with by hand.
\begin{corollary}
Let $n = 2^R + 2^{R-1}$ where $R \geq 2$. Let $\mu\vdash 2^{R-1}$ be an odd-dimensional partition of $m$ and $\lambda$ be a $2^R$-parent of $\mu$. Then
\[
\Od(f^{\lambda}) = (-1)^{s_2(h^\lambda_\mu) + \eta^\lambda_{\mu}}\Od(f^\mu).
\]
\end{corollary}

The proof of the above corollary is immediate as the sum of first two digits of $n$, $s_2(n)$, is 2. It follows from \cref{max-of-H-lambda} that $\mu$ is a hook partition of $2^{R-1}$. Equivalently, there exists a $b$ between 0 and $2^{R-1}-1$ such that $H(\mu) = \{2^{R-1}, b, b-1,\ldots, 1\}$. We first look at Type I $2^R$-parents which we denote by $\p_1(\mu)$ as in \cref{notn:parents}.
\begin{lemma}\label{lem:7}
Let $\mu\vdash 2^{R-1}$ be an odd-dimensional partition with $H(\mu) = \{2^{R-1},b, b-1,\ldots, 1\}$ for $R\geq 2$ and $b\geq 0$. Let $\p_1(\mu)$ denote the set of $2^R$-parents of $\mu$. We have 
\[
S_{\p_1(\mu)}(\mu) = \begin{cases}
2 ,& \text{if } \ell(\mu)\text{ is even}\\
1, & \text{if }\ell(\mu) \text{ is odd}.\\
\end{cases}
\]
\end{lemma}
\begin{proof}
For each $x\in H(\mu)$, define the partition $\la^x$ such that $H(\la^x)$ is obtained by replacing $x$ by $x+2^R$ in $H(\mu)$. We start by looking at $\nu:= \la^{2^{R-1}}$ for which $H(\nu)$ is obtained by replacing $2^{R-1}$ by $2^R + 2^{R-1}$ in $H(\mu)$. We deduce that $h^\nu_\mu = 2^R + 2^{R-1}$ and $s_2(h^\nu_\mu) = 2$. We compute $h^\nu_\mu - 2^{R-1} = 2^R$, $h^\nu_\mu - 2^R - 2^{R-1} = 0$, and $h^\nu_\mu+ 2^{R-1} = 2^{R+1}$.
Plugging these into \cref{prop:inversions}, we have

\[
\eta^\nu_\mu = N_\nu(2^R + 2^{R-1}) - \mathbb{I_{H(\nu)}}(2^R) + \mathbb{I_{H(\nu)}}(2^{R+1}) + \mathbb{I_{H(\nu)}}(0).
\]

As $2^{R-1}$ is the largest entry in $H(\mu)$, there cannot be any elements between $2^R + 2^{R-1}$ and $2^{R-1}$ in $H(\nu)$. Thus,
\[N_\nu(2^R + 2^{R-1}) = |\{y\in H(\nu) :  2^{R-1}< y < 2^{R}+2^{R-1}\}|=0.\] 
As $H(\mu)$ does not contain 0, $2^R$, or $2^{R+1}$, $H(\nu) = \{2^R + 2^{R-1}, b, \ldots, 1\}$ also does not contain any of these entries. The above formula for $\eta^\nu_\mu$ simplifies to $N_\nu(2^R + 2^{R-1})=0$. It follows that  $\Od(f^\nu) = \Od(f^\mu)$. 

Now we consider $\la^x$ for  $1\leq x \leq b$. As each $x$ is strictly smaller than $2^{R-1}$, $x + 2^R < 2^{R} + 2^{R-1}$ and $s_2(h^{\la^x}_\mu) = 1$. From \cref{prop:inversions}, it follows that
\begin{align*}
\eta^{\la^x}_\mu &= N_{\la^x}(2^R + x) - \mathbb{I}_{H(\la^x)}(x + 2^{R-1}) + \mathbb{I}_{H(\la^x)}(x + 2^R + 2^{R-1}) \\&\qquad\qquad+ \mathbb{I}_{H(\la^x)}(x - 2^{R-1}).
\end{align*}
Explicitly, $H(\la^x) = \{x + 2^R, 2^{R-1},b,\ldots, x+1, x-1, \ldots ,1\}$. There are $b+1-x$ entries in $H(\la^x)$ that lie strictly between $x+2^R$ and $x$, namely $x+1, x+2, \ldots, b$ and $2^{R-1}$.  From this, we deduce that $N_{\la^x}(2^R + x) = b+1-x$. By looking at the elements of $H(\mu)$, we can conclude that $H(\la^x)$ does not contain $ x - 2^{R-1}$, $x+2^{R-1}$ or $x + 2^{R} + 2^{R-1}$. Thus, $\eta^{\la^x}_\mu = N_{\mu}^{\la^x} = b + 1 - x$. 
For all elements $x\in H(\mu)$ not equal to $2^{R-1}$, $\Od(f^{\la^x})/\Od(f^\mu) = (-1)^{\eta^{\la^x}_\mu}+1$ where $\eta^{\la^x}_\mu = b-x + 1$ and the $+1$ appears in the exponent due to $s_2(h^{\la^x}_\mu) = 1$. When $x = 2^{R-1}$, $\Od(f^{\la^x})/\Od(f^\mu) = (-1)^{0+2}$ following from our previous computations.

Now we get (as in \cref{def:signed-sum}), 

\[
S_{\p_1(\mu)}({\mu}) = 1 + \sum\limits_{x=1}^{b} (-1)^{b + 1 - x + 1} = 1 + \sum\limits_{x=1}^{b} (-1)^{b - x} 
= 1 + \sum\limits_{y=0}^{b-1} (-1)^y
\]
where the last equality follows by setting $y = b-x$ in the sum. As we have elements $1,2,\ldots, b$ in $H(\mu)$, the value $b$ is one less than $\ell(\mu)$ and has the opposite parity to $\ell(\mu)$. This means that the upper bound $b-1$ has the same parity as that of $\ell(\mu)$. Thus, if the parity of $\ell(\mu)$ is even, $S_{\p_1(\mu)}(\mu) = 1+1 =2$ and if the parity of $\ell(\mu)$ is odd, $S_{\p_1(\mu)}(\mu) = 1+0 =1$.
\end{proof}
\begin{remark}
    We may interpret the above result in terms of hooks according to the conditions in \cref{rem:hook-interpretation}. We notice that $\mathbb{I}_{H(\la^x)}(h^{\la^x}_\mu - 2^{R}-2^{R-1}) = 0$ as hooks cannot be extended to non-rim hook subsets of size $2^R + 2^{R-1}$. This obstruction arises because such an extension of a rim hook requires all cells of $\la^x$ to be on the border of $\la^x$ and for them to not form a  hook shape, but that is not possible in our case. Also, as we are considering Type I parents, we find $\mathbb{I}_{H(\la^x)}(h^{\la^x}_\mu +2^{R-1}) = 0$, as the hooks cannot be extended upwards. Now, we notice that all hooks $\la^x/\mu$ can be decomposed into two $2^{R-1}$-rim hooks as the top row of $\la^x/\mu$ always contains at least $2^{R-1}+1$ boxes, which can be shown by size considerations. This makes $\mathbb{I}_{H(\la^x)}(h^{\la^x}_\mu -2^{R-1}) = 0$. In the following example, $\mu = (2,1,1)$:
    \[
    \sh(\la^1) = \y{2,1,1}*[*(gray)]{2+8}\quad \sh(\la^2) = \y{2,1,1}*[*(gray)]{2+6,1+2} \quad \sh(\la^3) = \y{2,1,1}*[*(gray)]{2+5,1+2,1+1}
    \]
\end{remark}
We now compute the signed-sum $S_{\p_2(\mu)}(\mu)$ for Type II $2^R$-parents of $\mu$. We focus on the values of $r$ for which $H(\mu)^{+r}$ contains $2^{R-1}$, $2^{R}$ or $2^R + 2^{R-1}$. We see in the following proof that the condition whether these numbers belong to $H(\mu)^{+r}$ leads to  the splitting of $[1,2^R]$ into six disjoint intervals.
\begin{lemma}\label{lem:8}
Let $\mu\vdash 2^{R-1}$ be an odd-dimensional partition of $2^{R-1}$ with $H(\mu) = \{2^{R-1},b, b-1,\ldots, 1\}$ for $R\geq 2$ and $b\geq 0$. Let $\p_2(\mu)$ denote the set of Type II $2^R$-parents of $\mu$, then 
\[
S_{\p_2(\mu)}({\mu}) = \begin{cases}
2 ,& \text{if }\ell(\mu)\text{ is even}\\
3, & \text{if }\ell(\mu) \text{ is odd}.\\
\end{cases}
\]
\end{lemma}
\begin{proof}
Denote by $\lambda_{[r]}$ the partition for which $H(\la_{[r]})$ is obtained by replacing 0 with $2^R$ in $ H(\mu)^{+r}$ when $2^R\not\in H(\mu)^{+r}$. $\p_2(\mu)$ is the set of such $\la_{[r]}$. For an invalid shift $s$ (\cref{notn:valid-values}), we let $\la_{[s]}$ be $\mu$. First, for all $\lambda\in\p_2(\mu)$, the affected hook-length (\cref{notn:aff}), $h^\lambda_\mu$, is  $2^R$, and so $s_2(h^\lambda_\mu) = 1$. This reduces the computation of the signed sum as follows:
\[
S_{\p_2(\mu)}({\mu}) = \sum\limits_{\lambda\in \p_2(\mu)} (-1)^{\eta^\lambda_\mu + 1} = (-1)\cdot\sum\limits_{\lambda\in \p_2(\mu)} (-1)^{\eta^\lambda_\mu}.
\]
As $2^R - (2^R+ 2^{R-1}) < 0$,   \cref{prop:inversions} tells us that $
\eta^\lambda_\mu = N_{\lambda}(2^R) - \mathbb{I}_{H(\lambda)}(2^{R-1}) + \mathbb{I}_{H(\lambda)}(2^R +2^{R-1})$ for any $\la\in \p_2(\mu)$.

In order to deduce whether $2^{R-1}$, $2^R$, or $2^{R}+2^{R-1}$ occur in $H(\la)$, we divide the set $\{1,2,\ldots, 2^R\}$ of values that $r$ can take into six intervals, of which some may be empty. As one may check, the following cases work for all values of $b\geq 0$.
\begin{enumerate}
\item Consider the case when {$1 \leq r\leq 2^{R-1} - b - 1$}. The largest entry in $H(\mu)^{+r}$ is $2^{R-1} + r$ which is at most $ 2^{R} - b -1$ and is thus strictly smaller than $2^R$. The second largest element in $H(\mu)^{+r}$ must be $b+r$ which is strictly smaller than $2^{R-1}$. Therefore, the maximum element of $H(\la_{[r]})$ is $2^R$ and $H(\la_{[r]})$ does not contain $2^{R-1}$. This forces $\mathbb{I}_{H(\lambda)}(2^{R-1})$ and $ \mathbb{I}_{H(\lambda)}(2^R +2^{R-1})$ to be zero. So, we have $\eta^{\lambda_{[r]}}_\mu = N_{\lambda_{[r]}}(2^R)$ which is the number of entries smaller than $2^R$. This is exactly one less than the number of entries in $H(\la)$. The length of $\la$ is $r$ more than the length of $\mu$ which gives us $\eta^{\lambda_{[r]}}_\mu = \ell(\mu) + r - 1$.

\item Consider {$2^{R-1} - b \leq r \leq 2^{R-1}-1$}. If $b = 0$, then the interval in this case is empty. If $b$ is positive, then $H(\mu)^{+r}$ contains $2^{R-1}$ which gives $\mathbb{I}_{H(\lambda)}(2^{R-1})=1$. On the other hand, as the largest value in $H(\mu)^{+r}$ is $b+r$, which is strictly smaller than $2^R + 2^{R-1}$, $\mathbb{I}_{H(\lambda)}(2^R + 2^{R-1})=0$. So $\eta^{\lambda_{[r]}}_\mu = N_{\lambda_{[r]}}(2^R) - 1$. We find $N_{\lambda_{[r]}}(2^R) \ell(\mu) + r-1$ as in part 1 and obtain $\eta^{\lambda_{[r]}}_\mu = \ell(\mu) + r - 2$.
\item When $r = 2^{R-1}$, we get $2^{R-1} + 2^{R-1} = 2^R\in H(\mu)^{+2^{R-1}}$, so $\lambda_{[2^{R-1}]} \not\in \p_2(\mu)$.

\item In the case when $2^{R-1} + 1 \leq r \leq 2^{R}-b-1$, $H(\lambda_{[r]})$ always contains $2^{R-1}$, so $\mathbb{I}_{H(\lambda)}(2^{R-1}) = 1$. Furthermore, the largest element of $H(\la_{[r]})$ is $2^{R-1} + r\leq 2^{R-1} + 2^R - b -1$ which is strictly smaller than $2^R + 2^{R-1}$, and thus $\mathbb{I}_{H(\lambda)}(2^R +2^{R-1}) = 0$. The number of elements in $H(\la_{[r]})$ are $\ell(\mu) + r$, out of which one element is $2^R$ and the other element is $2^{R-1} + 2^{R} - b - 1$, which is greater than $2^R$. This shows that the number of entries in $H(\la_{[r]})$ strictly smaller than $2^R$, that is $N_{\la_{[r]}}(2^R)$, is $\ell(\mu) + r - 2$. From this, we obtain $\eta^{\lambda_{[r]}}_\mu = (\ell(\mu)+r -2) - 1 = \ell(\mu) + r - 3$.

\item When $2^R - b \leq r \leq 2^{R}-1$, $2^{R}\in H(\mu)^{+r}$, so $\lambda_{[r]}\not\in \p_2(\mu)$.
\item Finally, when $r = 2^R$, the set $H(\mu)^{+2^R}$ contains both $2^{R-1}$ and $2^R + 2^{R-1}$. Furthermore, the set contains $1,2,\ldots, 2^{R}-1$, which means $N_{\la_{[2^R]}}(2^R) = 2^R - 1$.
So, $\eta_{\mu}^{\lambda_{[2^R]}} = N_{\lambda_{[2^R]}}(2^R) = (2^{R} - 1) -1 + 1 = 2^R -1$.
\end{enumerate}

We sum over the intervals in (1), (2), (4) and (6) to give
\begin{align*}
S_{\p_2(\mu)}({\mu}) &= (-1)\cdot\Bigg(\sum\limits_{r = 1}^{2^{R-1}-b-1} (-1)^{\ell(\mu)+r-1} + \sum\limits_{r = 2^{R-1}-b}^{2^{R-1}-1} (-1)^{\ell(\mu)+r-2} \\ 
&\qquad\qquad\quad+ \sum\limits_{r=2^{R-1} +1}^{2^{R}-b-1} (-1)^{\ell(\mu)+r - 3} + (-1)^{2^R - 1}\Bigg).
\end{align*}
We multiply the $-1$ inside the sum to find
\begin{align*}
S_{\p_2(\mu)}({\mu}) &= \sum\limits_{r = 1}^{2^{R-1}-b-1} (-1)^{\ell(\mu)+r} + \sum\limits_{r = 2^{R-1}-b}^{2^{R-1}-1} (-1)^{\ell(\mu)+r-1} \\ 
&\qquad\qquad\quad+ \sum\limits_{r=2^{R-1} +1}^{2^{R}-b-1} (-1)^{\ell(\mu)+r - 2} + (-1)^{2^R}.
\end{align*}
The last term is equal to $+1$. From the rest of the terms, we factor out $(-1)^{\ell(\mu)}$ to get

\[
 S_{\p_2(\mu)}({\mu}) = (-1)^{\ell(\mu)} \Bigg(\sum\limits_{r = 1}^{2^{R-1}-b-1} (-1)^{r} + \sum\limits_{r = 2^{R-1}-b}^{2^{R-1}-1} (-1)^{r-1} + \sum\limits_{r=2^{R-1} +1}^{2^{R}-b-1} (-1)^{r-2}\Bigg) + 1.
\]
The exponent $r-2$ in the third sum can be replaced with $r$ while keeping the same bounds. This yields
\[
 S_{\p_2(\mu)}({\mu}) = (-1)^{\ell(\mu)} \Bigg(\sum\limits_{r = 1}^{2^{R-1}-b-1} (-1)^{r} + \sum\limits_{r = 2^{R-1}-b}^{2^{R-1}-1} (-1)^{r-1} + \sum\limits_{r=2^{R-1} +1}^{2^{R}-b-1} (-1)^{r}\Bigg) + 1,
\]
and we notice that the third sum is the same as the first sum but the bounds are shifted up by $2^{R-1}$. Adding $2^{R-1}$ to the exponent of $(-1)^r$ does not change the sign of the output, so the first and third sums are equal. We now have 
\[
S_{\p_2(\mu)}({\mu})  = 1 + 2(-1)^{\ell(\mu)}\sum\limits_{r = 1}^{2^{R-1}-b - 1} (-1)^{r} + \sum\limits_{r = 2^{R-1}-b}^{2^{R-1}-1} (-1)^{r-1}.
\]
The $2^{R-1}$ in the bounds of the sum on the right hand side does not effect the value of the sum as $R\geq 2$. We reconsider this sum as ranging from $-b$ to $-1$, and the exponent of $-1$ ranges from $-(b+1)$ to $-2$. As $-1$ is its own inverse, we may write the sum as ranging over $x' = 2$ to $b+1$. A shift of 2 preserves the parity, and thus we may rewrite the sum as ranging over $x = 0$ to $x = b-1$.
We have 
\[
S_{\p_2(\mu)}({\mu})  = 1 + 2(-1)^{\ell(\mu)}\sum\limits_{r = 1}^{2^{R-1}-b - 1} (-1)^{r} + \sum\limits_{x=0}^{b-1} (-1)^{x}.
\]
The number of elements in $H(\mu)$ is one more than $b$, and so $b=\ell(\mu)-1$. This allows us to rewrite our sums in terms of $\ell(\mu)$ as
\[
S_{\p_2(\mu)}({\mu})  = 1 + 2(-1)^{\ell(\mu)}\sum\limits_{r = 1}^{2^{R-1}-\ell(\mu)} (-1)^{r} + \sum\limits_{x=0}^{\ell(\mu) - 2} (-1)^x.
\]
When $\ell(\mu)$ is even, we get $S_{\p_2(\mu)}({\mu}) = 1 + 2\cdot 1\cdot 0 + 1 = 2$. When $\ell(\mu)$ is odd, we get $S_{\p_2(\mu)}({\mu})  = 1 + 2\cdot (-1) \cdot (-1) + 0 = 3$.
\end{proof}
\begin{remark}
We interpret the intervals in terms of diagrams for $\gamma = \la/\mu$. The interval (1) corresponds to the case where $\gamma$ contains more than $2^{R-1}$ boxes in the top row of $\la$ and (2) corresponds to the case where $\gamma$ contains at most $2^{R-1}$, but still positive, number of boxes in the top row of $\la$. The interval (4) corresponds to $\gamma$ containing boxes weakly below the second row (if it exists) of $\la$, and (6) is when all boxes of $\gamma$ are contained in a column below the bottom row of $\la$.
\end{remark}
With all of these ingredients, we are ready to prove \cref{11thm}.
\begin{proof}[Proof of \cref{11thm}]
For a fixed $\mu\vdash 2^{R-1}$, we combine the results \[
S_{\p_1(\mu)}(\mu) = \begin{cases}
2 ,& \text{if } \ell(\mu)\text{ is even}\\
1, & \text{if }\ell(\mu) \text{ is odd}.\\
\end{cases}
\]
and \[
S_{\p_2(\mu)}({\mu}) = \begin{cases}
2 ,& \text{if }\ell(\mu)\text{ is even}\\
3, & \text{if }\ell(\mu) \text{ is odd}.\\
\end{cases}
\]
from \cref{lem:7} and \cref{lem:8} respectively, we see that $S_{\p(\mu)}(\mu) = 4$. Thus, for $n = 2^R + 2^{R-1}$, we find
\begin{align*}
\delta(n) &= \sum\limits_{\substack{\\\lambda\vdash n\\\lambda \text{ is odd}}} \Od(f^\lambda) = \sum\limits_{\substack{\mu\vdash 2^{R-1}\\\mu \text{ is odd}}} S_{\p(\mu)}({\mu})\Od(f^{\mu})\\
&= 4 \sum\limits_{\substack{\mu\vdash 2^{R-1}\\\mu \text{ is odd}}} \Od(f^{\mu})
= 4\delta(2^{R-1}).
\end{align*}
From \cref{lem:hookdims}, we know that $\delta(2^{R-1}) = 0$ for $R\geq 2$. This proves our theorem except when $n= 3$, which can be shown by explicitly computing $\delta(3)$ for the partitions $(1,1,1), (2,1)$ and $(3)$. Out of these $(1,1,1)$ and $(3)$ are odd-dimensional, and have dimensions equal to 1. So, $\delta(3) = 2$.
\end{proof}

\section{Tackling the case for general $n$}\label{sec:general-n}
For this section, we suppose that $n = 2^R + m$ with $2^{R-1}\leq m < 2^R$. If we could find a recursive expression for $\delta(n)$ for such an $n$, then we can compute the values of $\delta(k)$ for any natural number $k$ by using the results established in the previous sections. Such a recursive expression remains elusive. In this case, $s_2(n) = 2$, and so for any $2^R$-parent $\la$ of an odd-dimensional partition $\mu$ of $m$,
\[
\frac{\Od(f^\la)}{\Od(f^\mu)} = (-1)^{s_2(h^\la_\mu) + \eta^\la_\mu}. 
\]
In the previous cases, we saw that $S_{\p(\mu)}(\mu)$ only depended on the size of $\mu$ and not on the specific partition $\mu$. In the unsolved case for the general $n$, this does not always hold true. 
%For instance, let $\mu = (17,3,1)$ and $\nu = (12,3,3,2,1)$ be two partitions of $21 = [10101]_2$. Constructing their 32-parents, one can check that $S_{\p(\mu)}(\mu) = 8$ and $S_{\p(\nu)}(\nu) = 0$. 
Consider the following table for odd-dimensional partitions of 6 with $2^R = 8$. Each odd-dimensional partition of 6 has dimension 1 modulo 4, and so we find $\delta(14) = 0$.

\begin{table}[H]
    \centering
    \begin{tabular}{c|c|c}
        $\mu$ & $H(\mu)$ & $S_{\p(\mu)(\mu)}$ \\
        \hline
         (6)& \{6\}& 4\\
         (5,1)& \{6,1\}& $-4$\\
         (4,2)& \{5,2\}& $-4$\\
         (3,3)& \{4,3\}& 4\\
         (2,2,2) & \{4,3,2\} & 4\\
         (2,2,1,1) & \{5,4,2,1\} & $-4$\\
         (2,1,1,1,1) & \{6,4,3,2,1\} & $-4$\\
         (1,1,1,1,1,1) & \{6,5,4,3,2,1\} & 4.
    \end{tabular}
\end{table}

\begin{remark}
    It is not the case that $S_{\p(\mu)}(\mu)$ is always either $4$, $-4$ or 0. When $\mu = (31,15,4,3,1)\vdash 54$ with $2^R = 64$, we have $S_{\p(\mu)}(\mu) = 12$.
\end{remark}

Suppose $\mu$ is an odd-dimensional partition of $m$ and $\la$ is a Type I $2^R$-parent of $\mu$ obtained by replacing some $h\in H(\mu)$ by $h+2^R$. Then by \cref{prop:inversions},
\[
\eta^\la_\mu = N_\la(h+2^R) - \mathbb{I}_{H(\la)}(h + 2^{R-1}) + \mathbb{I}_{H(\la)}(h + 2^R + 2^{R-1})  + \mathbb{I}_{H(\la)}(h -2^{R-1}).
\]
When $h > 2^{R-1}$, 
$\eta^\la_\mu = N_\la(h+2^R) + \mathbb{I}_{H(\la)}(h -2^{R-1})$ and when $h < 2^{R-1}$, 
$\eta^\la_\mu = N_\la(h+2^R) - \mathbb{I}_{H(\la)}(h + 2^{R-1})$. In both these cases, the value of $\eta^\la_\mu$ is dictated by the presence of a pair of elements separated by $2^{R-1}$. Let $\nu$ be a Type II $2^R$-parent of $\mu$ obtained by replacing 0 by $2^R$ in some $r$-shift of $H(\mu)$. Then 
\[
\eta^\nu_\mu = N_\nu(2^R) - \mathbb{I}_{H(\nu)}(2^{R-1}) + \mathbb{I}_{H(\nu)}(2^R + 2^{R-1}).
\]
Suppose there exists an $h$ such that $h+r = 2^{R-1}\in H(\mu)$. As $\nu$ is a $2^R$-parent, then there must not exist an $h+2^{R-1}$ in $H(\mu)$.

These observations present a possible strategy to attack the general case. We categorize the elements of $H(\mu) = \{h_1', \ldots, h'_s, h_1, \ldots, h_u\}$ into two parts such that $h'_i \geq 2^{R-1}$ and $h_j < 2^{R-1}$. Under this segregation, we can pair up elements which differ by $2^{R-1}$, and it is by recording the contributions of these internally interacting pairs we can handle the general case.
In the case handled in \cref{sec:proofmain} all entries of the $\beta$-set were strictly smaller than $2^{R-1}$ and thus $s = 0$ and $u = \ell(\mu)$. On the other hand, for the case handled in \cref{sec:proof11}, we had exactly one entry equal to $2^{R-1}$ while the rest were smaller, so in our above notation $s = 1$ while $u = \ell(\mu) - 1$. In both these cases, there are no entries which are separated by $2^{R-1}$, which made our analysis simpler.

We end this section with one important observation. Assume $n = 2^{R} + m$ with $2^{R-1}\leq m < 2^R$. Let $\mu$ be an odd-dimensional partition of $m$, then $\mu$ contains a unique $2^{R-1}$-rim hook, which means that there exists some $1\leq t \leq s$ such that $h_t' - 2^{R-1}$ is not in $H(\mu)$. This implies that for all other $h_i'$ with $i\neq t$, there exists some $1\leq v_i \leq u$ such that $h_i' - 2^{R-1} = h_{v_i}$. So, we can remove a $2^{R-1}$-rim hook from $h'_t$ but not from other $h_i'$. In the case covered in \cref{sec:proof11}, we only had $t =u= 1$, which made it easier to analyze. We claim that it is due to the varying distribution of these internally interacting pairs $(h_i', h_{v_i})$ for $1\leq i \neq t\leq s$, that the values of $S_{\p(\mu)}$ vary for $\mu$.

\section{Acknowledgments}
I thank Steven Spallone for introducing me to this fascinating problem, and for his comments that helped streamline my mathematical thought and writing. I am grateful for the discussions I had with Jyotirmoy Ganguly which led to the insights that eventually became the proofs in this paper. My sincerest thanks go to Nick Loehr for his meticulous and consistent feedback. I also thank the anonymous reviewer for their detailed feedback on the manuscript that greatly improved its exposition.
\section{APPENDIX: Data for $\delta(n)$}\label{sec:problems}
We saw that $\delta(n)$ was a power of 4 for sparse numbers but this does not hold in general as can be seen in the table below. %\href{https://drive.google.com/file/d/1iaQz0ffXwzNc79LFg1_fDRkH-BLWTtL7/view?usp=drive_link}{data},
In general, $\delta(n)$ may not be even be a power of 2 as $\delta(118) = -384$. The run from 124 to 127 shows how egregious this deviation can get. We can show that $\delta(n) = a_1(n) - a_3(n)$ is even by noticing that for $n\geq 2$, $a(n) = a_1(n) + a_3(n)$ is even by McKay's result, and thus $a_1(n)$ and $a_3(n)$ have the same parity.

\[
\begin{tabular}{c|c}
$n$ & $\delta(n)$\\
\hline
 0 & 1 \\ 1 & 1 \\ 2 & 2 \\ 3 & 2 \\ 4 & 0 \\ 5 & 4\\ 6 & 8 \\ 7 & 0 \\ 8 & 0\\ 9 & 4 \\ 10 & 0\\ 11& 8 \\ 12 & 0 \\ 13 & 0 \\ 14 & 0 \\ 15 & 0 
\end{tabular}\qquad
\begin{tabular}{c|c}
    $n$ & $\delta(n)$ \\
    \hline
     16 & 0\\ 17 & 4 \\ 18 & 0 \\ 19 & 8 \\ 20 & 0 \\ 21 & 16 \\ 22 & 0 \\ 23 & 0 \\ 24 & 0 \\ 25 & 0 \\ 26 & 0 \\ 27 & 0 \\ 28 & 0 \\ 29 & 32 \\ 30 & $-64$ \\ 31 & 0
\end{tabular}\qquad
\begin{tabular}{c|c}
$n$ & $\delta(n)$\\
\hline
32 & 0 \\ 33 & 4 \\ 34 & 0 \\ 35 & 8 \\
36 & 0 \\ 37 & 16 \\ 38 & 0 \\ 39 & 0 \\ 40 & 0 \\ 41 & 16 \\ 42 & 0 \\ 43 & 32 \\ 44 & 0\\ 45 & 0\\ 46 & 0\\ 47 & 0
\end{tabular}\qquad
\begin{tabular}{c|c}
   $n$  &  $\delta(n)$ \\
     \hline
     48 & 0\\ 49 & 0\\ 50 & 0\\ 51 & 0\\ 52 & 0\\ 53 & 0\\ 54 & 0\\ 55 & 0\\ 56 & 0\\ 57 & 32\\ 58 & $-64$\\ 59 & 64\\60 & 64\\ 61 & $-128$\\ 62 & $-256$\\ 63 & 128
\end{tabular}
\]
\[
\begin{tabular}{c|c}
$n$ & $\delta(n)$\\
\hline
64 & 0\\ 65 & 4\\ 66 & 0\\ 67 & 8\\ 68 & 0\\ 69 & 16\\ 70  &  0\\
71  &  0\\
72  &  0\\73  &  16\\74  &  0\\75  &  32\\76  &  0\\77  &  0\\78  &  0\\79  &  0
\end{tabular}\qquad 
\begin{tabular}{c|c}
    $n$ & $\delta(n)$\\
     \hline
     80  &  0\\ 81  &  16\\82  &  0\\83  &  32\\84  &  0\\85  &  64\\86  &  0\\87  &  0\\88  &  0\\89  &  0\\90  &  0\\91  &  0\\92  &  0\\93  &  128\\94  &  0\\95  &  0
\end{tabular}\qquad
\begin{tabular}{c|c}
$n$ & $\delta(n)$\\
\hline
96  &  0\\97  &  0\\98  &  0\\99  &  0\\100  &  0\\101  &  0\\102  &  0\\103  &  0\\104  &  0\\105  &  0\\106  &  0\\107  &  0\\
108  &  0\\
109  &  0\\110  &  0\\ 111  &  0
\end{tabular}\qquad
\begin{tabular}{c|c}
    $n$ & $\delta(n)$\\
     \hline
     112  &  0\\ 113  &  32\\114  &  $-64$\\115  &  64\\116  &  $-64$\\117  &  64\\118  &  $-384$\\119  &  128\\120  &  256\\121  &  64\\122  &  $-256$\\123  &  $-256$\\124  &  768\\125  &  640\\126  &  168\\127  &  256 
\end{tabular}
\]

\bigskip
{\bf \large Contact Information}\\
%Write email addresses, affiliations, and ORCID of the authors. Then Comment them out using % for double-blind review
%ORCID icon used below can be downloaded from 
%https://ajcombinatorics.org/ORCIDiD_icon128x128.png

\begin{tabular}{lcl}
%%Author 1	
Aditya Khanna	&& Virginia Tech\\
adityakhanna$@$vt.edu && 225 Stanger St, Blacksburg, VA 24060, USA\\
	&& 
    \url{https://orcid.org/0009-0002-1933-5313}\\
				&& \\
                \end{tabular}
%%%%%%%%%%%%%%%%%%%%%%%%%%%%%%%%%%%				

\end{document}